\documentclass[12pt]{article}
\usepackage{a4wide,color}
\usepackage[centertags,leqno]{amsmath}
\usepackage{amssymb}
\usepackage[mathscr]{eucal}
\usepackage{epsf}

\newtheorem{theorem}{Theorem}[section]

\newtheorem{example}[theorem]{Example}
\newtheorem{remark}[theorem]{Remark}
\newtheorem{cor}[theorem]{Corollary}

\newtheorem{definition}[theorem]{Definition}

\newtheorem{hypo}{Hypothesis}

\newcommand{\qed}{\rule{2mm}{2mm}}

\newcommand{\eqdef}{\stackrel{{\mathrm {def}}}{=}}
\newcommand{\dotV}{\overset{\bullet}{V}}
\newcommand{\eps}{\varepsilon}
\renewcommand{\colon}{:\,}
\newcommand{\proof}{{\em Proof. }}

\newcommand{\RR}{\mathbb{R}}

\renewcommand{\SS}{\mathbb{S}}

\numberwithin{equation}{section}

\begin{document}

\title{\bf
       A $p(x)$\--Laplacian Extension of the D\'{\i}az\--Saa\\
       Inequality and Some Applications\\
\vspace{1.5cm}
}

\author{%
\vspace{0.1cm}
        Peter {\sc Tak\'a\v{c}}\\
        Institut f\"ur Mathematik,
        Universit\"at Rostock\\
        Ulmenstra{\ss}e~69, Haus~3\\
        D-18055 Rostock, Germany\\
{\tt peter.takac@uni-rostock.de}\\
\and
        Jacques {\sc Giacomoni}\\
        LMAP (UMR 5142)\\
        Universit\'e de Pau et des Pays de l'Adour\\
        Avenue de l'Universit\'e,
        F-64013 Pau cedex, France\\
{\tt jacques.giacomoni@univ-pau.fr}\\
\vspace{0.1cm}
}

\maketitle

\vspace{0.5cm}

\newpage
\noindent
\begin{abstract}
\baselineskip=14pt
The main result of this work is a new extension of
the well\--known inequality by D\'{\i}az and Saa which, in our case,
involves an anisotropic operator, such as the $p(x)$\--Laplacian,
$\Delta_{p(x)} u\equiv \mathrm{div} (|\nabla u|^{p(x) - 2} \nabla u)$.
Our present extension of this inequality enables us to establish
several new results on the uniqueness of solutions
and comparison principles
for some anisotropic quasilinear elliptic equations.
Our proofs take advantage of certain convexity properties of
the energy functional associated with the $p(x)$\--Laplacian.
\end{abstract}

\vfill
\par\vspace*{0.5cm}
\noindent
\begin{tabular}{ll}
{\bf Running head:}
& D\'{\i}az and Saa Inequality for $\Delta_{p(x)}$\\
\end{tabular}

\par\vspace*{0.5cm}
\noindent
\begin{tabular}{ll}
{\bf Keywords:}
& $p(x)$-Laplacian;\\
& quasilinear Dirichlet problem with variable exponents;\\
& ray\--strictly convex energy functional;\\
& uniqueness and comparison principles\\
\end{tabular}

\par\vspace*{0.5cm}
\noindent
\begin{tabular}{lll}
{\bf 2000 Mathematics Subject Classification:}
& Primary   & 35J62, 35J92;\\
& Secondary & 35B09, 35A02 \\
\end{tabular}


\baselineskip=16pt

\newpage
\section{Introduction}
\label{s:Intro}

This work is concerned with an extension of
a well\--known inequality due to
{\sc J.~I.\ D{\'\i}az} and {\sc J.~E.\ Saa} \cite{Diaz-Saa}
to certain quasilinear elliptic operators that are pointwise
$p(x)$\--homogeneous, but anisotropic, in general, such as
the $p(x)$\--Laplacian
\begin{math}
  \Delta_{p(x)} u\equiv\hfil\break
  \mathrm{div} (|\nabla u|^{p(x) - 2} \nabla u)
\end{math}
with a variable exponent $p(x)\in (1,\infty)$.
Such operators have been studied extensively, e.g., in
{\sc L.\ Diening}, {\sc P.\ Harjulehto}, {\sc P.\ H\"ast\"o}, and
{\sc M.\ R\r{u}\v{z}i\v{c}ka} \cite{DienHHR}, and in
{\sc V.\ R\u{a}dulescu} and {\sc D.\ Repov\v{s}} \cite{Radu-Rep}.
Interesting applications to a model of {\it electrorheological fluids\/}
are presented in \cite[{\S}14.4, pp.\ 470--481]{DienHHR}
and the monograph by
{\sc M.\ R\r{u}\v{z}i\v{c}ka} \cite{Ruzicka}.
However, to our best knowledge,
the {\sc D{\'\i}az} and {\sc Saa} inequality \cite{Diaz-Saa}
still has not been extended from the original case of
a constant exponent
$p(x)\equiv p = \mathrm{const}\in (1,\infty)$
to a {\it variable\/} exponent $p(x)$.
This inequality turns out to be equivalent with the convexity of
a $p(x)$\--power type energy functional, as suggested in
{\sc H.\ Br\'ezis} and {\sc L.\ Oswald} \cite{Brezis-Osw}
for $p(x)\equiv p=2$, and generalized in
{\sc J.\ Fleckinger}, {\sc J.\ Hern\'andez}, {\sc P.\ Tak\'a\v{c}},
and {\sc F.\ de~Th\'elin} \cite{FleckHTT}
to any constant $p(x)\equiv p\in (1,\infty)$.
In applications to quasilinear elliptic operators
(with $p$ constant, $1 < p < \infty$),
this equivalence played a decisive role in the works by
{\sc P.\ Girg} and {\sc P.\ Tak\'a\v{c}}
\cite[{\S}4.1, pp.\ 289--292]{Girg-Takac}
and
{\sc P.\ Tak\'a\v{c}}, {\sc L.\ Tello}, and {\sc M.\ Ulm}
\cite[Lemma 2.4, p.~79]{TakacUT}.

To be more specific, the functional in question,
$\mathcal{W}\colon W\to \RR_+$, is defined by
\begin{equation}
\label{def:W(v)}
  \mathcal{W}(v)\equiv \mathcal{W}_{p(x),r}(v)\eqdef
  \int_{\Omega} \frac{r}{p(x)}\,
  \left\vert \nabla (|v(x)|^{1/r}) \right\vert^{p(x)} \,\mathrm{d}x
\end{equation}
for every function $v\in L^{p(x)/r}(\Omega)$ such that
$|v|^{1/r}\in W_0^{1,p(x)}(\Omega)$;
the set of all such functions $v\colon \Omega\to \RR$ is denoted by
$W\equiv W_{p(x),r}$.
Here, we assume that $\Omega\subset \RR^N$ is a bounded domain in
$\RR^N$ ($N\geq 1$) whose boundary is a compact manifold
(with smoothness to be specified later if $N\geq 2$),
$r\in [1,\infty)$ is a given constant, and
$p\in L^{\infty}(\Omega)$ is an essentially bounded function satisfying
$p(x) > 1$ and $p(x)\geq r$ for almost all $x\in \Omega$
(whose smoothness will be specified later, as well).
We will show in the next section
(Section~\ref{s:Main})
that this functional is {\it convex\/} on the cone
\begin{equation}
\label{def:V_cone}
  \dotV\eqdef \{ v\colon \Omega\to (0,\infty)\colon v\in W\}
  \subset W
\end{equation}
of all positive functions $v\in W$.
The convexity of the restriction
$\mathcal{W}\colon \dotV\to \RR_+$ to $\dotV$
is well\--known to be equivalent with the monotonicity of its
(set\--valued) subdifferential $\partial\mathcal{W}(v)$ at $v\in \dotV$
that is a nonempty set only for certain elements $v\in \dotV$
which might not be easy to determine; cf.\
{\sc P.\ Girg} and {\sc P.\ Tak\'a\v{c}}
\cite[{\S}4.1, pp.\ 289--292]{Girg-Takac}.
To avoid this problem, we restrict ourselves only to
certain directional derivatives of $\mathcal{W}$
which exist in the classical sense.

We consider two functions $v_1, v_2\in \dotV$ such that
$v_1/v_2$, $v_2/v_1\in L^{\infty}(\Omega)$.
Consequently, also
$v\eqdef (1-\theta) v_1 + \theta v_2\in \dotV$
is valid for all $\theta\in (-\delta, 1+\delta)$,
where $\delta\in (0,1)$ is small enough.
The function
\begin{equation*}
  \theta\mapsto \Phi(\theta)\eqdef \mathcal{W}(v)
  = \mathcal{W}\left( (1-\theta) v_1 + \theta v_2\right)
    \colon (-\delta, 1+\delta)\to \RR_+
\end{equation*}
is convex and differentiable with the derivative
\begin{equation}
\label{e:dW/d_theta}
  \Phi'(\theta) =
  \int_{\Omega}
  \left\vert \nabla (|v(x)|^{1/r}) \right\vert^{p(x) - 2}
  \nabla (|v(x)|^{1/r})\cdot
  \nabla \genfrac{(}{)}{}0{v_2-v_1}{ v^{ 1 - \frac{1}{r} } }
  \,\mathrm{d}x \,;
\end{equation}
see Theorem~\ref{thm-ray-convex} below.
The monotonicity of the derivative
\begin{math}
  \theta\mapsto \Phi'(\theta)\colon (-\delta, 1+\delta)\to \RR
\end{math}
yields $\Phi'(1) - \Phi'(0)\geq 0$, i.e.,
\begin{equation}
\label{ineq:dW/d_theta}
\begin{aligned}
& \left\langle
{}- \frac{ \Delta_{p(x)} (v_1(x)^{1/r}) }{ v_1(x)^{(r-1)/r} }
  + \frac{ \Delta_{p(x)} (v_2(x)^{1/r}) }{ v_2(x)^{(r-1)/r} }
  \,,\; v_1 - v_2
  \right\rangle
\\
& = \int_{\Omega}
  \left(
{}- \frac{ \Delta_{p(x)} (v_1(x)^{1/r}) }{ v_1(x)^{(r-1)/r} }
  + \frac{ \Delta_{p(x)} (v_2(x)^{1/r}) }{ v_2(x)^{(r-1)/r} }
  \right)
  (v_1 - v_2) \,\mathrm{d}x\geq 0 \,,
\end{aligned}
\end{equation}
provided the integration by parts in eq.~\eqref{e:dW/d_theta}
can be justified. 
In this case we may substitute
$w_i = v_i^{1/r} > 0$ in $W_0^{1,p(x)}(\Omega)$; $i=1,2$,
to derive the following extension of
the {\it\bfseries D{\'\i}az and Saa inequality\/}
(Theorem~\ref{thm-Diaz-Saa} below):
\begin{equation}
\label{ineq:Diaz-Saa}
  \int_{\Omega}
  \left(
{}- \frac{ \Delta_{p(x)} w_1(x) }{ w_1(x)^{r-1} }
  + \frac{ \Delta_{p(x)} w_2(x) }{ w_2(x)^{r-1} }
  \right)
  (w_1^r - w_2^r) \,\mathrm{d}x\geq 0
\end{equation}
for all pairs $w_1, w_2\in W_0^{1,p(x)}(\Omega)$, such that
$w_1 > 0$, $w_2 > 0$ a.e.\ in $\Omega$ and both
$w_1/w_2$, $w_2/w_1\in L^{\infty}(\Omega)$.
The special case
$p(x)\equiv r = \mathrm{const}\in (1,\infty)$
yields the classical {\sc D{\'\i}az} and {\sc Saa} inequality
established in \cite{Diaz-Saa}.

To verify ineq.~\eqref{ineq:Diaz-Saa},
in Section~\ref{s:Main} below (Theorem~\ref{thm-Diaz-Saa})
we slightly modify the method used in
\cite{FleckHTT, Girg-Takac, TakacUT}.
Our proof of ineq.~\eqref{ineq:Diaz-Saa}
is based on the convexity of the restriction of
the functional $\mathcal{W}$ to the cone $\dotV\subset W$.
Finally, in Section~\ref{s:Appl}
we present a few applications of our main results
to some nonlinear boundary value problems with
the Dirichlet $p(x)$\--Laplacian $\Delta_{p(x)}$ and
the power\--type nonlinearity $|u(x)|^{q(x)-2} u(x)$, where
the following (uniform) {\it ``subhomogeneity'' condition\/}
is assumed:
\begin{equation}
\label{e:subhomog}
  1 < q(x)\leq r\leq p(x)
  \quad\mbox{ for almost every }\, x\in \Omega \,,
\end{equation}
with a suitable (uniform separation) constant $r\in (1,\infty)$.
This condition is related to abstract subhomogeneity conditions
introduced in the well\--known monograph by
{\sc M.~A.\ Krasnosel'ski\u{\i}} and {\sc P.~P. Zabre\u{\i}ko}
\cite{Krasno-Zabr}
in several different abstract settings in ordered Banach spaces.

\section{Main Results and Their Proofs}
\label{s:Main}

It is easy to see that the set $\dotV$ defined in eq.~\eqref{def:V_cone}
is a {\em convex cone\/}, i.e.,
{\rm (i)}$\;$
$\lambda\in (0,\infty)$, $f\in \dotV$ $\,\Rightarrow\,$
$\lambda f\in \dotV$; and
{\rm (ii)}$\;$
$f,g\in \dotV$ $\,\Rightarrow\,$ $f+g\in \dotV$.

\begin{definition}\label{def-ray-convex}\nopagebreak
\begingroup\rm
A functional\/ $\mathcal{W}\colon \dotV\to \RR$
will be called {\em\bfseries ray\--strictly convex\/}
({\em\bfseries strictly convex\/}, respectively)
if it satisfies
\begin{align}
\label{e:ray-strictly}
  \mathcal{W}\left( (1-\theta) v_1 + \theta v_2 \right)
  \leq (1-\theta)\cdot \mathcal{W}(v_1)
      +   \theta \cdot \mathcal{W}(v_2)
\\
\nonumber
    \quad\mbox{ for all $v_1, v_2\in \dotV$ and for all }
                \theta\in (0,1) \,,
\end{align}
where the inequality is strict ($<$) unless
$v_2 / v_1\equiv \mathrm{const} > 0$ is a constant
(for $\mathcal{W}$ ``strictly convex'' always strict if\/ $v_1\neq v_2$).
\endgroup
\end{definition}
\par\vskip 10pt

We assume that $\Omega\subset \RR^N$ is
either a bounded open interval in $\RR^1$ ($N=1$)
or else a bounded domain in $\RR^N$ ($N\geq 2$)
whose boundary $\partial\Omega$ is a compact manifold of class
$C^{1,\alpha}$ for some $\alpha\in (0,1)$.
Additional hypotheses on the smoothness of the boundary $\partial\Omega$
(such as interior sphere condition at $\partial\Omega$)
will be added later in the applications
(Section~\ref{s:Appl}).

For the sake of simplicity, we assume that
$r\in [1,\infty)$ is a given constant and
$p\colon \Omega\to (1,\infty)$ is a continuous function, such that
\begin{equation}
\label{e:p-<p<p+}
  1 < p_{-}\eqdef \inf_{\Omega} p(x)\leq
      p_{+}\eqdef \sup_{\Omega} p(x) < \infty
  \quad\mbox{ and }\quad 1\leq r\leq p_{-} \,.
\end{equation}

We assume that the function $A$ of
$(x,\xi)\in \Omega\times \RR^N$ extends to
a continuous and nonnegative function
$A\colon \overline{\Omega}\times \RR^N\to \RR_+$,
and it verifies the following hypothesis:
For every fixed $x\in \Omega$, the function
$A(x,\,\cdot\,)\colon \RR^N\to \RR_+$
is {\em positively $p(x)$\--homo\-ge\-neous\/}, i.e.,
\begin{equation}
\label{p-hom:A}
  A(x,t\xi) = |t|^{p(x)}\, A(x,\xi) \quad\mbox{ for all }\,
    t\in \RR \,,\ \xi\in\RR^N \,.
\end{equation}

It is evident from eq.~\eqref{p-hom:A} that
\begin{equation*}
  A(x,\xi) =
  A\left( x,\, \frac{\xi}{|\xi|} \right)\cdot |\xi|^{p(x)} \,,
    \quad\mbox{ for all }\, \xi\in \RR^N\setminus \{\mathbf{0}\} \,,
\end{equation*}
determines the growth of $A(x,\xi)$ with respect to
$\xi\in \RR^N\setminus \{\mathbf{0}\}$,
for any fixed $x\in \Omega$.
Let
$\SS^{N-1}\eqdef \{ \xi\in \RR^N\colon |\xi| = 1\}$
denote the unit sphere in $\RR^N$ centered at the origin
$\mathbf{0}\in \RR^N$.
We remark that the ``coefficient''
\begin{math}
  A\left( x,\, \frac{\xi}{|\xi|} \right)
\end{math}
in the last equation above
is bounded from above by a positive constant,
thanks to
$A\colon \overline{\Omega}\times \SS^{N-1}\to \RR_+$
being continuous on the compact set
\begin{math}
  \overline{\Omega}\times \SS^{N-1}\subset \RR^N\times \RR^N \,.
\end{math}

The simpliest example of $A$ is, of course,
$A(x,\xi) = |\xi|^{p(x)}$ for $(x,\xi)\in \Omega\times \RR^N$,
in which case $A(x,\xi) = 1$ for all
$(x,\xi)\in \Omega\times \SS^{N-1}$.
This case leads to the functional $\mathcal{W}$ defined in
eq.~\eqref{def:W(v)}.

The next theorem is our main result on the functional
$\mathcal{W}_{A}\colon W\to \RR_+$ defined by
\begin{align}
\label{def:W_A(v)}
  \mathcal{W}_A(v)
& \equiv \mathcal{W}_{A,p(x),r}(v)\eqdef
  \int_{\Omega} \frac{r}{p(x)}\,
  A\left( x,\, \nabla (|v(x)|^{1/r}) \right) \,\mathrm{d}x
\\
\nonumber
& {}
  = \int_{\Omega} \frac{r}{p(x)}\,
    A \left( x,\, \frac{ \nabla (|v|^{1/r}) }%
                  { \left\vert \nabla (|v|^{1/r}) \right\vert }
      \right)\cdot
    \left\vert \nabla (|v|^{1/r}) \right\vert^{p(x)} \,\mathrm{d}x
\end{align}
for every function $v\in W$; see eq.~\eqref{def:W(v)}
in the Introduction (Section~\ref {s:Intro}).

\begin{theorem}\label{thm-ray-convex}
{\rm (Convexity)}$\;$
Let\/ $r\in [1,\infty)$ and\/
$p\colon \Omega\to (1,\infty)$ satisfy \eqref{e:p-<p<p+}.
Assume that\/
\begin{math}
  A\colon \overline{\Omega}\times \RR^N\to \RR_+
\end{math}
is continuous and satisfies
the $p(x)$\--homogeneity hypothesis \eqref{p-hom:A}.
In addition, assume that the function
\begin{equation}
\label{e:A-conv}
  \xi\mapsto \mathfrak{N}(x,\xi)\eqdef A(x,\xi)^{r/p(x)}
  \colon \RR^N\to \RR_+
\end{equation}
is strictly convex for every $x\in \Omega$.
Then the restriction of the functional\/
$\mathcal{W}_A\colon W\to \RR_+$ to the convex cone $\dotV$ is
{\bf ray\--strictly convex\/} on~$\dotV$.

Furthermore, if\/ $p(x)\not\equiv r$ in $\Omega$, i.e., if\/
$r = p_{-}\equiv p(x)\equiv p_{+}$ does not hold in~$\Omega$,
then $\mathcal{W}_A$ is even
{\bf strictly convex\/} on~$\dotV$.
\end{theorem}
\par\vskip 10pt

\begin{remark}\label{rem-ray-convex}\nopagebreak
\begingroup\rm
{\rm (i)}$\;$
In the classical setting with $p(x)\equiv p\in (1,\infty)$
being a constant and $r=1$
(cf.\ {\sc P.\ Tak\'a\v{c}} \cite{Takac-5}),
\begin{math}
  \mathfrak{N}(x,\xi)\equiv \mathfrak{N}(\xi)
  = |\xi| = \|\xi\|_{\ell_2}
\end{math}
($\xi\in \RR^N$)
is the standard Euclidean norm in $\RR^N$.
Hence, the functional
$\mathcal{W}_{A}\colon W\to \RR_+$ defined in
\eqref{def:W(v)} and \eqref{def:W_A(v)} takes the form
\begin{align*}
  \mathcal{W}_A(v)\equiv \mathcal{W}_{A,p,1}(v)
  = \frac{1}{p}\int_{\Omega} |\nabla v(x)|^p \,\mathrm{d}x
  = \frac{1}{p}\, \| v\|_{ W_0^{1,p}(\Omega) }^p
\end{align*}
for every $v\in W = W_0^{1,p}(\Omega)$, thanks to
$|\nabla |v|| = |\nabla v|$ a.e.\ in $\Omega$.

{\rm (ii)}$\;$
In fact, Part~{\rm (i)} was the motivation for expressing the function
\begin{math}
  \xi\mapsto A(x,\xi) =\hfil\break
  \mathfrak{N}(x,\xi)^{p(x)/r} \colon \RR^N\to \RR_+
\end{math}
as a power ($p(x)/r\geq 1$) of the (strictly convex) function
\begin{math}
  \xi\mapsto \mathfrak{N}(x,\xi) \colon \RR^N\to \RR_+
\end{math}
that may be taken to be a strictly convex norm on $\RR^N$
depending on $x\in \Omega$; cf.\
{\sc P.\ Tak\'a\v{c}}, {\sc L.\ Tello}, and {\sc M.\ Ulm}
\cite[Remark 2.1, p.~78]{TakacUT}.

{\rm (iii)}$\;$
We note that
(in Theorem~\ref{thm-ray-convex} above)
the function
\begin{math}
  \xi\mapsto A(x,\xi) =\hfil\break
  \mathfrak{N}(x,\xi)^{p(x)/r} \colon \RR^N\to \RR_+
\end{math}
is strictly convex for each fixed $x\in \Omega$,
thanks to the power function
$t\mapsto t^{p(x)/r}\colon \RR_+\to \RR_+$
being strictly monotone increasing and convex.
Consequently,
$A(x,\xi) > A(x,\mathbf{0}) = 0$ for all $x\in \Omega$ and
$\xi\in \RR^N\setminus \{\mathbf{0}\}$, and
$A\colon \overline{\Omega}\times \SS^{N-1}\to \RR_+$
is bounded below and above on the compact set
\begin{math}
  \overline{\Omega}\times \SS^{N-1}\subset \RR^N\times \RR^N
\end{math}
by some positive constants; hence, the ``coefficient''
\begin{math}
  A\left( x,\, \frac{\xi}{|\xi|} \right)
\end{math}
in the integrand in eq.~\eqref{def:W_A(v)}, if
$\xi = \nabla (|v|^{1/r})\neq \mathbf{0}$,
is bounded from below and above by some positive constants
$c_1, c_2\in (0,\infty)$,
\begin{equation*}
  0 < c_1\leq A\left( x,\, \genfrac{}{}{}1{\xi}{|\xi|} \right)
         \leq c_2 < \infty \,.
\end{equation*}
This shows that also the ratio of the functionals in
\eqref{def:W_A(v)} and \eqref{def:W(v)}
is bounded from below and above by the same positive constants as above,
i.e.,
\begin{equation*}
  c_1\cdot \mathcal{W}(v)\leq \mathcal{W}_A(v)\leq
  c_2\cdot \mathcal{W}(v)
  \quad\mbox{ for every }\, v\in \dotV \,.
\end{equation*}
\endgroup
\end{remark}
\par\vskip 10pt

{\it Proof of\/} Theorem~\ref{thm-ray-convex}.
Recalling {\rm Definition~\ref{def-ray-convex}},
let us consider any $v_1, v_2\in \dotV$ and $\theta\in (0,1)$.
Let us denote
$v = (1-\theta) v_1 + \theta v_2$; hence, $v\in \dotV$.
We obtain easily
\begin{align*}
    \nabla (v_i(x)^{1/r})
& = \frac{v_i^{1/r}}{r}\, \frac{\nabla v_i}{v_i}
    \quad\mbox{ for }\, i=1,2 \,, \quad\mbox{ and }
\\
    \nabla (v(x)^{1/r})
& = \frac{1}{r}\,
    \frac{ (1-\theta) \nabla v_1 + \theta\nabla v_2 }%
         { [ (1-\theta) v_1 + \theta v_2 ]^{1 - (1/r)} }
\\
  = \frac{v^{1/r}}{r}\,
    \frac{ (1-\theta) \nabla v_1 + \theta\nabla v_2 }{v}
& = \frac{v^{1/r}}{r}\,
    \left[ (1-\theta)\, \frac{v_1}{v}\cdot \frac{ \nabla v_1 }{v_1}
             + \theta\, \frac{v_2}{v}\cdot \frac{ \nabla v_2 }{v_2}
   \right] \,,
\end{align*}
with the convex combination of positive coefficients
$(1-\theta)\, \frac{v_1}{v}$ and $\theta\, \frac{v_2}{v}$,
\begin{equation*}
  (1-\theta)\, \frac{v_1}{v} + \theta\, \frac{v_2}{v} = 1 \,.
\end{equation*}
Now let $x\in \Omega$ be fixed.
Since $\xi\mapsto \mathfrak{N}(x,\xi)$
is strictly convex, by our hypothesis,
we may apply the identities from above to conclude that
\begin{equation}
\label{e:v<v_1+v_2}
\begin{aligned}
& \mathfrak{N}
  \left(x, (1-\theta)\, \frac{v_1}{v}\cdot \frac{ \nabla v_1 }{v_1}
            + \theta\,  \frac{v_2}{v}\cdot \frac{ \nabla v_2 }{v_2}
  \right)
\\
& \leq (1-\theta)\, \frac{v_1}{v}\cdot
  \mathfrak{N}
  \left(x, \frac{ \nabla v_1 }{v_1} \right)
        + \theta\, \frac{v_2}{v}\cdot
  \mathfrak{N}
  \left(x, \frac{ \nabla v_2 }{v_2} \right) \,.
\end{aligned}
\end{equation}
The equality holds if and only if
\begin{equation}
\label{e:v_1=v_2}
    \frac{ \nabla v_1(x) }{v_1(x)}
  = \frac{ \nabla v_2(x) }{v_2(x)} \,,
    \quad\mbox{ which is equivalent to }\quad
  \nabla\genfrac{(}{)}{}0{v_2(x)}{v_1(x)} = 0 \,.
\end{equation}
Notice that the homogeneity conditions
\eqref{p-hom:A} and \eqref{e:A-conv} yield
\begin{equation}
\label{p-hom:N}
  \mathfrak{N}(x,t\xi) = |t|^r\, \mathfrak{N}(x,\xi)
    \quad\mbox{ for all }\, t\in \RR \,,\ \xi\in\RR^N \,.
\end{equation}
Consequently, we multiply ineq.~\eqref{e:v<v_1+v_2} by
$v / r^r$ to obtain the following equivalent inequality,
\begin{equation}
\label{eq:v<v_1+v_2}
\begin{aligned}
& \mathfrak{N}
    \left( x,\nabla (v(x)^{1/r})\right)
  = \frac{v}{r^r}\cdot
  \mathfrak{N}
  \left( x, (1-\theta)\, \frac{v_1}{v}\cdot \frac{ \nabla v_1 }{v_1}
              + \theta\, \frac{v_2}{v}\cdot \frac{ \nabla v_2 }{v_2}
  \right)
\\
& \leq (1-\theta)\, \frac{v_1}{r^r}\cdot
  \mathfrak{N}
  \left( x, \frac{ \nabla v_1 }{v_1} \right)
      + \theta\,   \frac{v_2}{r^r}\cdot
  \mathfrak{N}
  \left( x, \frac{ \nabla v_2 }{v_2} \right)
\\
& = (1-\theta)\cdot \mathfrak{N}\left( x, \nabla (v_1(x)^{1/r})\right)
      + \theta\cdot \mathfrak{N}\left( x, \nabla (v_2(x)^{1/r})\right) \,.
\end{aligned}
\end{equation}

Finally, by {\rm Remark~\ref{rem-ray-convex}},
we conclude that ineq.~\eqref{eq:v<v_1+v_2} entails
\begin{equation}
\label{eq:A:v<v_1+v_2}
\begin{aligned}
&   A\left( x, \nabla (v(x)^{1/r})\right)
  \leq
    (1-\theta)\cdot A\left( x, \nabla (v_1(x)^{1/r})\right)
      + \theta\cdot A\left( x, \nabla (v_2(x)^{1/r})\right) \,.
\end{aligned}
\end{equation}

We multiply the last inequality, \eqref{eq:A:v<v_1+v_2},
by $r/p(x)$, then integrate the product over $\Omega$ to derive
the convexity of the restriction of
the functional $\mathcal{W}_A$ to the convex cone $\dotV\subset W$.

To derive that $\mathcal{W}_A$
is even ray\--strictly convex on $\dotV$, let us consider any pair
$v_1, v_2\in \dotV$ with $v_1\not\equiv v_2$ in~$\Omega$.
We observe that the equality in the convexity inequality
\eqref{e:ray-strictly} forces both conditions,
\eqref{e:v_1=v_2} and $p(x)/r = 1$, to hold simultaneously
at almost every point $x\in \Omega$.
These conditions are then equivalent with
$v_2/v_1\equiv \mathrm{const}$ ($\neq 1$) in $\Omega$ and
$p(x)\equiv r$ in $\Omega$.
Thus, if $p(x)\not\equiv r$ in $\Omega$,
then $\mathcal{W}_A$ is even strictly convex on~$\dotV$.
\qed
\par\vskip 10pt

Our second theorem is concerned with the extension of
the {\sc D{\'\i}az} and {\sc Saa} {\it inequality\/}
as formulated in ineq.~\eqref{ineq:Diaz-Saa}.
Here, we need to assume a more specific form of the function
\begin{math}
  A\colon \overline{\Omega}\times \RR^N\to \RR_+ .
\end{math}
Besides the homogeneity hypothesis \eqref{e:A-conv},
we assume that $A$ and its partial gradient
\begin{math}
  \partial_\xi A\equiv
  \left( \frac{\partial A}{\partial\xi_i} \right)_{i=1}^N
\end{math}
with respect to $\xi\in \RR^N$ satisfy
the following structural hypothesis, upon the substitution
\begin{math}
  \mathbf{a}(x,\xi) \eqdef \frac{1}{p(x)}\, \partial_\xi A(x,\xi)
\end{math}
with
\begin{math}
  a_i = \frac{1}{p(x)}\, \frac{\partial A}{\partial\xi_i} :
\end{math}
%

\begin{hypo}\nopagebreak
\begingroup\rm
({\bf A})
$\;$
Given any fixed $x\in \Omega$, the function
$A(x,\,\cdot\,)\colon \RR^N\to \RR_+$
verifies the {\em positive $p(x)$\--homo\-ge\-nei\-ty\/} hypothesis
\eqref{p-hom:A}.
Furthermore, we assume that
\begin{math}
  A\in\hfil\break
  C(\overline{\Omega}\times \RR^N)\cap C^1(\Omega\times \RR^N)
\end{math}
and its partial gradient
$\partial_\xi A\colon \Omega\times \RR^N\to \RR^N$
satisfies
\begin{math}
  \frac{1}{p}\, \frac{\partial A}{\partial\xi_i} = a_i
  \in C^1( \Omega\times (\RR^N \setminus \{ 0\}) )
\end{math}
for all $i=1,2,\dots,N$, together with
the following {\em ellipticity\/} and {\em growth conditions\/}:
There exist some constants
$\gamma,\Gamma\in (0,\infty)$ such that
%
\begin{eqnarray}
  \sum_{i,j=1}^N
  \frac{\partial a_i}{\partial\xi_j} (x,\xi)
    \cdot \eta_i\eta_j
&\geq &
  \gamma \cdot |\xi|^{p(x)-2} \cdot |\eta|^2 ,
\label{h:ellipt}
\\
  \sum_{i,j=1}^N
    \left|
  \frac{\partial a_i}{\partial\xi_j} (x,\xi)
    \right|
&\leq &
  \Gamma \cdot |\xi|^{p(x)-2} ,
\label{h:growth-xi}
\end{eqnarray}
for all $x\in \Omega$, all $\xi\in \RR^N\setminus \{ 0\}$,
and all $\eta\in \RR^N$.
\endgroup
\end{hypo}
\par\vskip 10pt

Owing to the homogeneity hypothesis \eqref{e:A-conv},
it suffices to assume that the inequalities in
\eqref{h:ellipt} and \eqref{h:growth-xi}
hold for all $\xi\in \SS^{N-1}$ only.

\begin{theorem}\label{thm-Diaz-Saa}
\begingroup\rm
(The {\sc D{\'\i}az} and {\sc Saa} inequality.)$\;$
\endgroup
Let\/ $r\in [1,\infty)$ and\/
$p\colon \Omega\to (1,\infty)$ satisfy \eqref{e:p-<p<p+}.
Assume that\/
\begin{math}
  A\colon \overline{\Omega}\times \RR^N\to \RR_+
\end{math}
satisfies {\rm Hypothesis ({\bf A})} and, in addition, the function
\begin{math}
  \xi\mapsto \mathfrak{N}(x,\xi) = A(x,\xi)^{r/p(x)}
  \colon \RR^N\to \RR_+
\end{math}
is strictly convex for every $x\in \Omega$.
Then the following inequality\/
\begin{equation}
\label{in_A:Diaz-Saa}
  \int_{\Omega}
  \left(
{}- \frac{ \mathrm{div}\, \mathbf{a}(x, \nabla w_1(x)) }{ w_1(x)^{r-1} }
  + \frac{ \mathrm{div}\, \mathbf{a}(x, \nabla w_2(x)) }{ w_2(x)^{r-1} }
  \right)
  (w_1^r - w_2^r) \,\mathrm{d}x\geq 0
\end{equation}
holds (in the sense of distributions)
for all pairs $w_1, w_2\in W_0^{1,p(x)}(\Omega)$, such that
$w_1 > 0$, $w_2 > 0$ a.e.\ in $\Omega$ and both
$w_1/w_2$, $w_2/w_1\in L^{\infty}(\Omega)$.
Moreover, if the equality ($=$) in \eqref{in_A:Diaz-Saa} occurs,
then we have the following two statements:
\begin{itemize}
\item[{\rm (a)}]
$\;$
$w_2/w_1\equiv \mathrm{const} > 0$ in $\Omega$.
\item[{\rm (b)}]
$\;$
If also $p(x)\not\equiv r$ in $\Omega$, then even
$w_1\equiv w_2$ holds in $\Omega$.
\end{itemize}
\end{theorem}
\par\vskip 10pt

\begin{remark}\label{rem-Diaz-Saa}\nopagebreak
\begingroup\rm
The distributional inequality~\eqref{in_A:Diaz-Saa}
has to be interpreted in the following way:
\begin{equation}
\label{in_A:Diaz-Saa:int}
\begin{aligned}
& \int_{\Omega} \mathbf{a}(x, \nabla w_1(x))\cdot \nabla
  \left( w_1 - \frac{w_2^r}{w_1^{r-1}} \right) \,\mathrm{d}x
\\
& \geq 
  \int_{\Omega} \mathbf{a}(x, \nabla w_2(x))\cdot \nabla
  \left( \frac{w_1^r}{w_2^{r-1}} - w_2 \right) \,\mathrm{d}x
\end{aligned}
\end{equation}
for all pairs $w_1, w_2\in W_0^{1,p(x)}(\Omega)$, such that
$w_1 > 0$, $w_2 > 0$ a.e.\ in $\Omega$ and both
$w_1/w_2$, $w_2/w_1\in L^{\infty}(\Omega)$.
Both integrals above are defined as Lebesgue integrals,
thanks to the inequalities in
\eqref{h:ellipt} and \eqref{h:growth-xi}
combined with the following standard identities,
\begin{align*}
    \nabla \genfrac{(}{)}{}0{ w_2^r }{ w_1^{r-1} }
  = r \genfrac{(}{)}{}0{w_2}{w_1}^{r-1} \nabla w_2
  - (r-1) \genfrac{(}{)}{}0{w_2}{w_1}^r \nabla w_1 \,,
\\
    \nabla \genfrac{(}{)}{}0{ w_1^r }{ w_2^{r-1} }
  = r \genfrac{(}{)}{}0{w_1}{w_2}^{r-1} \nabla w_1
  - (r-1) \genfrac{(}{)}{}0{w_1}{w_2}^r \nabla w_2 \,,
\end{align*}
where
$w_1/w_2$, $w_2/w_1\in L^{\infty}(\Omega)$
and all gradients belong to $L^{p(x)}(\Omega)$, whence also
\begin{math}
  { w_2^r }/{ w_1^{r-1} } ,\,\hfil\break
  { w_1^r }/{ w_2^{r-1} } \in W_0^{1,p(x)}(\Omega) .
\end{math}
\endgroup
\end{remark}
\par\vskip 10pt

{\it Proof of\/} Theorem~\ref{thm-Diaz-Saa}.
Recalling {\rm Definition~\ref{def-ray-convex}},
let us consider any pair
$w_1, w_2\in W_0^{1,p(x)}(\Omega)$, such that
$w_1 > 0$, $w_2 > 0$ a.e.\ in $\Omega$ and both
$w_1/w_2$, $w_2/w_1\in L^{\infty}(\Omega)$.
Consequently, there is a number $\delta\in (0,1)$,
sufficiently small, such that
\begin{equation*}
  v\eqdef (1-\theta) w_1^r + \theta w_2^r\in \dotV
    \quad\mbox{ and }\quad
  v^{1/r}\in W_0^{1,p(x)}(\Omega)
    \;\mbox{ for all }\, \theta\in (-\delta, 1+\delta) \,.
\end{equation*}
The function
\begin{equation*}
  \theta\mapsto \Phi(\theta)\eqdef \mathcal{W}(v)
  = \mathcal{W}_A\left( (1-\theta) w_1^r + \theta w_2^r\right)
    \colon (-\delta, 1+\delta)\to \RR_+
\end{equation*}
is convex and differentiable with the derivative
\begin{equation}
\label{e:dW_A/d_theta}
  \Phi'(\theta) =
  \int_{\Omega} \mathbf{a}(x, \nabla (v(x)^{1/r}))\cdot
  \nabla \genfrac{(}{)}{}0{ w_2^r - w_1^r }{ v^{ 1 - \frac{1}{r} } }
  \,\mathrm{d}x \,.
\end{equation}
To provide a rigorous proof of the convexity claim,
one has to consider two arbitrary points
$\theta_1, \theta_2\in \RR$, such that
$-\delta < \theta_1 < \theta_2 < 1+\delta$,
and all their convex combinations
$\theta = (1-t)\theta_1 + t\theta_2\in (-\delta, 1+\delta)$
with $t\in [0,1]$.
For $0\leq \theta_1 < \theta_2\leq 1$
the convexity is known, by Theorem~\ref{thm-ray-convex}.
However, if at least one of the following inequalities holds,
$-\delta < \theta_1 < 0$ and/or
$1 < \theta_2 < 1+\delta$,
the convexity inequality
\begin{math}
  \Phi(\theta)\leq (1-t)\Phi(\theta_1) + t\Phi(\theta_2)
\end{math}
still remains to be verified.
Of course, the number $\delta > 0$ needs to be taken small enough.
We leave this easy exercise to the reader.

The monotonicity of the derivative
\begin{math}
  \theta\mapsto \Phi'(\theta)\colon (-\delta, 1+\delta)\to \RR
\end{math}
yields $\Phi'(0)\leq \Phi'(1)$, which is equivalent with
ineq.~\eqref{in_A:Diaz-Saa:int}, thanks to
$v =  w_1^r$ if $\theta = 0$, and
$v =  w_2^r$ if $\theta = 1$.
It is now easy to see that ineq.~\eqref{in_A:Diaz-Saa}
is a distributional interpretation of \eqref{in_A:Diaz-Saa:int}
after integration by parts.

Finally, let us assume that
the equality ($=$) in \eqref{in_A:Diaz-Saa} is valid.
This forces $\Phi'(0) = \Phi'(1)$ above; hence,
$\Phi'(\theta) = \Phi'(0)$ for all $\theta\in [0,1]$,
by the monotonicity of $\Phi'\colon [0,1]\to \RR$.
It follows that $\Phi\colon [0,1]\to \RR$ must be linear, i.e.,
$\Phi(\theta) = (1-\theta) \Phi(0) + \theta \Phi(1)\in \RR$
for all $\theta\in [0,1]$.
Recalling our definition of $\Phi$ above and
Theorem~\ref{thm-ray-convex},
we conclude that $w_2/w_1\equiv \mathrm{const} > 0$ in $\Omega$.
This proves statement {\rm (a)}.

To verify statement {\rm (b)}, suppose that the constant above
$w_2/w_1\equiv \mathrm{const} \neq 1$ in $\Omega$.
Then the equality in both inequalities,
\eqref{eq:v<v_1+v_2} and \eqref{eq:A:v<v_1+v_2},
is possible only if $p(x)\equiv r$ in $\Omega$.
Statement {\rm (b)} follows.
\qed
\par\vskip 10pt

Our third (and last) theorem is a {\em weak comparison principle\/}
for positive solutions $u\in W_0^{1,p(x)}(\Omega)$
of the following (uniformly) {\it ``subhomogeneous''\/}
Dirichlet boundary value problem:
\begin{equation}
\label{e:BVP}
\left\{
\begin{aligned}
{}- \mathrm{div}\, \mathbf{a}(x, \nabla u(x))
& = f(x)\, u(x)^{r-1}
    \quad\mbox{ for }\, x\in \Omega \,;\qquad
  u > 0 \;\mbox{ a.e. in }\, \Omega \,;
\\
  u&= 0 \quad\mbox{ for }\, x\in \partial\Omega \,.
\end{aligned}
\right.
\end{equation}
Here, $f\in L^{\infty}(\Omega)$ is a given nonnegative function,
$f\geq 0$ a.e.\ in $\Omega$.

\begin{theorem}\label{thm-weak_comp}
\begingroup\rm
(Weak comparison principle.)$\;$
\endgroup
Let all\/ $r\in [1,\infty)$,
$p\colon \Omega\to (1,\infty)$,
\begin{math}
  A\colon \overline{\Omega}\times \RR^N\to \RR_+ ,
\end{math}
and the function
\begin{math}
  \xi\mapsto \mathfrak{N}(x,\xi) = A(x,\xi)^{r/p(x)}
  \colon \RR^N\to \RR_+
\end{math}
satisfy the same hypotheses as in
{\rm Theorem~\ref{thm-Diaz-Saa}} above.
In addition, assume that\/
$p(x)\not\equiv r$ in $\Omega$, i.e.,
$p(x) > r$ on a subset of\/ $\Omega$ with positive Lebesgue measure.

Finally, let\/
$u_i\in W_0^{1,p(x)}(\Omega)$
be a positive solution of the Dirichlet boundary value problem
\eqref{e:BVP} (in the sense of distributions)
with $f = f_i\in L^{\infty}(\Omega)$ for $i=1,2$, respectively,
where $0\leq f_1\leq f_2$ a.e.\ in $\Omega$.
If\/
$u_1/u_2$, $u_2/u_1\in L^{\infty}(\Omega)$
then we have also\/
$u_1\leq u_2$ a.e.\ in $\Omega$.
\end{theorem}
\par\vskip 10pt

\begin{remark}\label{rem-weak_comp}\nopagebreak
\begingroup\rm
Conditions
$u_1/u_2$, $u_2/u_1\in L^{\infty}(\Omega)$
imposed on the solutions $u_1$ and $u_2$ are easy to verify for
$f_1\not\equiv 0$ (hence, also $f_2\not\equiv 0$) in $\Omega$,
by the regularity results in
{\sc X.\ Fan} and {\sc D.\ Zhao} \cite[Theorem 4.1, p.~312]{Fan-Zhao}
and
{\sc X.-L.\ Fan} \cite[Theorem 1.2, p.~400]{Fan-JDE}
combined with the Hopf boundary point lemma from
{\sc Q.\ Zhang} \cite[Theorems 1.1 and 1.2, p.~26]{Zhang};
see our proof of Theorem~\ref{thm-problem_1} below.
\endgroup
\end{remark}
\par\vskip 10pt

It will be obvious from our proof of Theorem~\ref{thm-weak_comp}
below that the following simple generalization of this theorem
to weak sub- and super\-solutions is a direct consequence of the proof.
(We leave the details concerning only the last two inequalities
 of the proof to an interested reader.)

\begin{theorem}\label{thm-sub-super}
\begingroup\rm
(Weak comparison principle for sub- and super\-solutions.)$\;$
\endgroup
Let all\/ $r\in [1,\infty)$,
$p\colon \Omega\to (1,\infty)$,
\begin{math}
  A\colon \overline{\Omega}\times \RR^N\to \RR_+ ,
\end{math}
and the function
\begin{math}
  \xi\mapsto \mathfrak{N}(x,\xi) = A(x,\xi)^{r/p(x)}
  \colon \RR^N\to \RR_+
\end{math}
satisfy the same hypotheses as in
{\rm Theorem~\ref{thm-Diaz-Saa}} above.
In addition, assume that\/
$p(x)\not\equiv r$ in $\Omega$, i.e.,
$p(x) > r$ on a subset of\/ $\Omega$ with positive Lebesgue measure.

Finally, let\/
$u_i\in W_0^{1,p(x)}(\Omega)$ ($i=1,2$)
be a pair of positive functions satisfying\/
$u_1/u_2$, $u_2/u_1\in L^{\infty}(\Omega)$
together with the following inequalities (in the sense of distributions)
with $f_i\in L^{\infty}(\Omega)$ for $i=1,2$, respectively,
where $0\leq f_1\leq f_2$ a.e.\ in $\Omega$:
\begin{align}
\label{e:BVP_1}
\left\{
\begin{aligned}
{}- \mathrm{div}\, \mathbf{a}(x, \nabla u_1(x))
& \leq f_1(x)\, u_1(x)^{r-1}
    \quad\mbox{ for }\, x\in \Omega \,;\quad
  u_1 > 0 \;\mbox{ a.e. in }\, \Omega \,;
\\
  u_1&= 0 \quad\mbox{ for }\, x\in \partial\Omega \,.
\end{aligned}
\right.
\\
\label{e:BVP_2}
\left\{
\begin{aligned}
{}- \mathrm{div}\, \mathbf{a}(x, \nabla u_2(x))
& \geq f_2(x)\, u_2(x)^{r-1}
    \quad\mbox{ for }\, x\in \Omega \,;\quad
  u_2 > 0 \;\mbox{ a.e. in }\, \Omega \,;
\\
  u_2&= 0 \quad\mbox{ for }\, x\in \partial\Omega \,.
\end{aligned}
\right.
\end{align}
Then also\/
$u_1\leq u_2$ a.e.\ in $\Omega$ holds.
\end{theorem}
\par\vskip 10pt

We quote a well\--known fact from the theory of distributions that
any nonnegative distribution in $\mathcal{D}'(\Omega)$
may be identified with a nonnegative Radon measure on $\Omega$.
This result shows that the left\--hand side of both inequalities
\eqref{e:BVP_1} and \eqref{e:BVP_2}
must be a Radon measure on $\Omega$.

\par\vskip 10pt
{\it Proof of\/} Theorem~\ref{thm-weak_comp}.
We proceed in analogy with the proof of
{\rm Theorem~\ref{thm-Diaz-Saa}} above.
We set $w_i = u_i$; $i=1,2$, and define
\begin{equation*}
  v\equiv v(\theta)\eqdef u_2^r + \theta (u_1^r - u_2^r)^{+}
    \;\mbox{ for all }\, \theta\in (-\delta, 1+\delta) \,,
\end{equation*}
where $\delta\in (0,1)$ is a sufficiently small number, such that
$v\in \dotV$ for every $\theta\in (-\delta, 1+\delta)$.
As usual, the symbol
$\xi^{+} = \max\{ \xi, 0\}\geq 0$
stands for the positive part of a real number $\xi\in \RR$.
Hence, we have also
$v^{1/r}\in W_0^{1,p(x)}(\Omega)$.
Notice that
\begin{equation*}
  v = \left\{
  \begin{array}{cl}
  u_2^r
  &\quad\mbox{ if }\, u_1\leq u_2 \,,
\\
  u_2^r + \theta (u_1^r - u_2^r)
  &\quad\mbox{ if }\, u_1 > u_2 \,.
  \end{array}
\right.
\end{equation*}
On the contrary to our claim
$u_1\leq u_2$ a.e.\ in $\Omega$, let us assume that
$(u_1^r - u_2^r)^{+} > 0$ holds on a subset
\begin{math}
  \Omega_{+} = \{ x\in \Omega\colon u_1(x) > u_2(x) \}\subset \Omega
\end{math}
of positive Lebesgue measure.

By Theorem~\ref{thm-ray-convex}, thanks to our hypothesis
$p(x)\not\equiv r$ in $\Omega$, the function
\begin{equation*}
  \theta\mapsto \Phi(\theta)\eqdef \mathcal{W}(v)
  = \mathcal{W}_A\left(  u_2^r + \theta (u_1^r - u_2^r)^{+} \right)
    \colon (-\delta, 1+\delta)\to \RR_+
\end{equation*}
is strictly convex and differentiable with the derivative
\begin{equation}
\label{eq:dW_A/d_theta}
  \Phi'(\theta) =
  \int_{\Omega} \mathbf{a}(x, \nabla (v(x)^{1/r}))\cdot \nabla
  \genfrac{(}{)}{}0{ (u_1^r - u_2^r)^{+} }{ v^{ 1 - \frac{1}{r} } }
  \,\mathrm{d}x \,.
\end{equation}
The strict convexity of
\begin{math}
  \Phi\colon (-\delta, 1+\delta)\to \RR
\end{math}
and the monotonicity of its derivative $\Phi'$ yield
$\Phi'(0) < \Phi'(1)$, which is equivalent with
\begin{equation*}
\begin{aligned}
& \int_{\Omega} \mathbf{a}(x, \nabla u_2(x))\cdot \nabla
  \genfrac{(}{)}{}0{ (u_1^r - u_2^r)^{+} }{ u_2^{r-1} }
  \,\mathrm{d}x
  <
\\
& \int_{\Omega}
  \mathbf{a}
  \left( x, \nabla\left[ ( u_2^r + (u_1^r - u_2^r)^{+} )^{1/r} \right]
  \right)
  \cdot \nabla
  \genfrac{(}{)}{}0{ (u_1^r - u_2^r)^{+} }%
                   { ( u_2^r + (u_1^r - u_2^r)^{+} )^{ 1 - \frac{1}{r} } }
  \,\mathrm{d}x \,.
\end{aligned}
\end{equation*}
By {\rm Remark~\ref{rem-Diaz-Saa}},
the last inequality has the following distributional interpretation,
\begin{equation}
\label{ineq:Phi'(0,1)}
\begin{aligned}
& {}
  - \int_{\Omega}
    \frac{ \mathrm{div}\; \mathbf{a}(x, \nabla u_2(x)) }{ u_2(x)^{r-1} }
  \, (u_1^r - u_2^r)^{+} \,\mathrm{d}x
  <
\\
& {}
  - \int_{\Omega}
    \frac{ \mathrm{div}\;
  \mathbf{a}
  \left( x, \nabla\left[ ( u_2^r + (u_1^r - u_2^r)^{+} )^{1/r} \right]
  \right)
  }{ ( u_2^r + (u_1^r - u_2^r)^{+} )^{ 1 - \frac{1}{r} } }
  \, (u_1^r - u_2^r)^{+} \,\mathrm{d}x \,.
\end{aligned}
\end{equation}
As it is well\--known from the theory of Sobolev spaces of type
$W_0^{1,p(x)}(\Omega)$, both integrands above vanish almost everywhere
in the Lebesgue measurable set
\begin{math}
  \Omega_{-} = \{ x\in \Omega\colon
  u_1(x)\leq \hfil\break u_2(x) \}\subset \Omega .
\end{math}
Consequently, ineq.~\eqref{ineq:Phi'(0,1)} reads
\begin{equation}
\label{in_+:Phi'(0,1)}
  \int_{\Omega_{+}}
  \left(
{}- \frac{ \mathrm{div}\, \mathbf{a}(x, \nabla u_2(x)) }{ u_2(x)^{r-1} }
  + \frac{ \mathrm{div}\, \mathbf{a}(x, \nabla u_1(x)) }{ u_1(x)^{r-1} }
  \right)
  (u_1^r - u_2^r) \,\mathrm{d}x < 0 \,.
\end{equation}
By our hypotheses, we have
\begin{equation*}
{}- \frac{ \mathrm{div}\, \mathbf{a}(x, \nabla u_2(x)) }{ u_2(x)^{r-1} }
  + \frac{ \mathrm{div}\, \mathbf{a}(x, \nabla u_1(x)) }{ u_1(x)^{r-1} }
  = f_2(x) - f_1(x)\geq 0 \quad\mbox{ for a.e. }\, x\in \Omega \,.
\end{equation*}
Since also $u_1^r - u_2^r > 0$ a.e.\ in $\Omega_{+}$,
ineq.~\eqref{in_+:Phi'(0,1)} leads to a contradiction.

Thus, we have proved that the set $\Omega_{+}$ must have
Lebesgue measure equal to zero.
\qed
\par\vskip 10pt

\section{Applications to Differential Equations}
\label{s:Appl}

In this section we give two applications of
Theorems \ref{thm-ray-convex} and~\ref{thm-Diaz-Saa}.
Throughout this section we impose the following hypotheses
on $\Omega$ and $p(x)$:

\begin{hypo}\nopagebreak
\begingroup\rm
($\boldsymbol{\Omega}$)
$\;$
If $N=1$ then $\Omega$ is a bounded open interval in $\RR^1$.
If $N\geq 2$ then $\Omega$ is a bounded domain in $\RR^N$
whose boundary $\partial\Omega$ is a compact manifold of class
$C^{1,\alpha}$ for some $\alpha\in (0,1)$, and
$\Omega$ satisfies also the interior sphere condition
at every point of $\partial\Omega$.
\endgroup
\end{hypo}
\par\vskip 10pt

It is clear that for $N\geq 2$,
Hypothesis {\rm ($\boldsymbol{\Omega}$)} is satisfied if, for instance,
$\Omega\subset \RR^N$ is a bounded domain with $C^2$ boundary.
We write
$\overline{\Omega} = \Omega\cup \partial\Omega$
for the closure of $\Omega$ in $\RR^N$.

\begin{hypo}\nopagebreak
\begingroup\rm
($\mathbf{p}$)
$\;$
We assume that
$p\colon \overline{\Omega}\to (1,\infty)$
is $\alpha_1$-H\"older\--continuous, i.e.,
$p\in C^{0,\alpha_1}(\overline{\Omega})$ for some $\alpha_1\in (0,1)$,
and $p$ satisfies \eqref{e:p-<p<p+} with a given constant
$r\in [1,\infty)$, i.e.,
\begin{equation*}
  1 < p_{-}\eqdef \inf_{\Omega} p(x)\leq
      p_{+}\eqdef \sup_{\Omega} p(x) < \infty
  \quad\mbox{ and }\quad 1\leq r\leq p_{-} \,.
\end{equation*}
\endgroup
\end{hypo}
\par\vskip 10pt

Our first application is the following nonlinear
Dirichlet boundary value problem taken from
{\sc L.\ Diening}, {\sc P.\ Harjulehto}, {\sc P.\ H\"ast\"o}, and
{\sc M.\ R\r{u}\v{z}i\v{c}ka}
\cite[Eq.\ (13.3.2), p.~418]{DienHHR},
\begin{equation}
\label{e:problem_1}
\left\{
\begin{alignedat}{2}
  {}- \Delta_{p(x)} u
& {}= f(x,u)
&&  \quad\mbox{ in }\, \Omega \,;
\\
  u & {}= 0
&&  \quad\mbox{ on }\, \partial\Omega \,,
    \quad u > 0 \mbox{ in }\, \Omega \,.
\end{alignedat}
\right.
\end{equation}

\par\vskip 10pt
%
We impose the following hypotheses on the function $f$:
\begin{itemize}
\item[{\bf (f1)}]
$\;$
$f\colon \overline{\Omega}\times \RR_+\to \RR_+$
is a nonnegative continuous function such that $f(x,0) = 0$ for all
$x\in \Omega$.
\item[{\bf (f2)}]
$\;$
The function
\begin{math}
  s\,\longmapsto\, f(x,s) / s^{r-1} \colon (0,\infty)\to \RR_+
\end{math}
is strictly monotone decreasing for every $x\in \Omega$.
\item[{\bf (f3)}]
$\;$
The following two limits are uniform with respect to $x\in \Omega$:
\begin{equation*}
  \frac{ f(x,s) }{ s^{r-1} } \,\longrightarrow\, +\infty
  \quad\mbox{ as }\, s\to 0+\,
  \quad\mbox{ and }\quad
  \frac{ f(x,s) }{ s^{r-1} } \,\longrightarrow\, 0
  \quad\mbox{ as }\, s\to +\infty \,.
\end{equation*}
Equivalently, we require
\begin{equation*}
\begin{alignedat}{2}
& \frac{1}{ s^{r-1} }\cdot\; \inf_{x\in \Omega} f(x,s)
  \,\longrightarrow\, +\infty
&&\quad\mbox{ as }\, s\to 0+\,
  \quad\mbox{ and }\quad
\\
& \frac{1}{ s^{r-1} }\cdot\; \sup_{x\in \Omega} f(x,s)
  \,\longrightarrow\, 0
&&\quad\mbox{ as }\, s\to +\infty \,.
\end{alignedat}
\end{equation*}
\end{itemize}
%
\par\vskip 10pt

A typical example of the function $f$ satisfying
all {\rm Hypotheses\/} {\bf (f1)} -- {\bf (f3)}, with
$f(x,s) = h(x)\, s^{q(x)-1}$ for $x\in \overline{\Omega}$
and $s\in \RR_+$, is given below in Example~\ref{exam-appl_1}.
Here, $h\in C(\overline{\Omega})$ is a positive function and
$q\in C(\overline{\Omega})$ satisfies
\begin{math}
  1\leq q(x)\leq q_{+}\eqdef \sup_{\Omega} q(x) < r = p_-
\end{math}
for every $x\in \overline{\Omega}$.
In fact, we may choose any number $r\in (q_+,p_-]$
while requiring $q_+ < p_-$.
As a consequence, in this example we must have
$1\leq q_+ < r\leq p_-$ whence $r > 1$.

We remark that Hypothesis~{\bf (f3)} implies
the following asymptotic behavior of the function
$s\mapsto f(x,s)\colon (0,\infty)\to \RR_+$
as $s\to 0+$:
Given any $\eps > 0$, there is a constant $s_{\eps}\in (0,\infty)$
such that
\begin{equation}
\label{lim:s=0}
  f(x,s)\geq \frac{1}{\eps}\, s^{r-1}
  \quad\mbox{ holds for all }\,
  (x,s)\in \overline{\Omega}\times [0,s_{\eps}] \,.
\end{equation}
In contrast, Hypotheses {\bf (f1)} and {\bf (f3)}
limit the asymptotic behavior of $f(x,\,\cdot\,)$
as $s\to +\infty$ as follows:
Given any $\eps > 0$, there is a constant $C_{\eps}\in (0,\infty)$
such that
\begin{equation}
\label{lim:s=infty}
  0\leq f(x,s)\leq \eps\, s^{r-1} + C_{\eps}
  \quad\mbox{ holds for all }\,
  (x,s)\in \overline{\Omega}\times \RR_+ \,.
\end{equation}

We define the notion of a {\it nonnegative\/} weak solution to
problem \eqref{e:problem_1} as follows:

\begin{definition}\label{def-weak_sol}\nopagebreak
\begingroup\rm
A nonnegative function
$u\in W^{1,p(x)}_0(\Omega)\cap L^\infty(\Omega)$
is called a {\em\bfseries nonnegative weak solution\/}
of problem \eqref{e:problem_1} if,
for every test function $\phi\in W^{1,p(x)}_0(\Omega)$,
the following equation holds,
\begin{equation}
\label{def:problem_1}
  \int_{\Omega}
    \vert\nabla u(x)\vert^{p(x)-2}\nabla u(x)\cdot \nabla\phi(x)
   \,\mathrm{d}x 
  = \int_{\Omega} f(x,u(x))\, \phi(x) \,\mathrm{d}x \,.
\end{equation}

If\/ $u$ satisfies also $u > 0$ throughout $\Omega$, we call $u$
a {\em\bfseries positive weak solution\/}.
\endgroup
\end{definition}
\par\vskip 10pt

Problem~\eqref{e:problem_1}
has already been treated in
{\sc X.-L.\ Fan} and {\sc Q.-H.\ Zhang} \cite{Fan-Zhang-1}
where the existence of a weak solution in $W^{1,p(x)}_0(\Omega)$
is proved; see also
{\sc L.\ Diening}, {\sc P.\ Harjulehto}, {\sc P.\ H\"ast\"o}, and
{\sc M.\ R\r{u}\v{z}i\v{c}ka}
\cite[Theorem 13.3.3, p.~418]{DienHHR}.
Of course, the trivial solution $u\equiv 0$ in $\Omega$
is a nonnegative weak solution to problem \eqref{e:problem_1}.

The following theorem describes the solvability of
the boundary value problem~\eqref{e:problem_1}
for positive weak solutions.

\begin{theorem}\label{thm-problem_1}
Under the {\rm Hypotheses\/}
{\rm ($\boldsymbol{\Omega}$)}, {\rm ($\mathbf{p}$)}, and\/
{\bf (f1)} -- {\bf (f3)},
problem~\eqref{e:problem_1} possesses
a unique nonnegative and nontrivial weak solution
$u\in W^{1,p(x)}_0(\Omega)\cap L^\infty(\Omega)$.
This solution belongs to the class
$C^{1,\beta}(\overline{\Omega})$, for some $\beta\in (0, \alpha)$,
and satisfies also the Hopf maximum principle,
\begin{equation}
\label{Hopf:problem_1}
  u(x) > 0 \;\mbox{ for all }\, x\in \Omega
  \quad\mbox{ and }\quad
  \frac{\partial u}{\partial\boldsymbol{\nu}} (x) < 0
    \;\mbox{ for all }\, x\in \partial\Omega \,.
\end{equation}
Of course, $u = 0$ on the boundary $\partial\Omega$.
Hence, $u$ is also a positive weak solution.
\end{theorem}
\par\vskip 10pt

As usual, the symbol $\boldsymbol{\nu}(x)\in \RR^N$ stands for
the {\em unit outward normal\/} to the boundary $\partial\Omega$
at the point $x\in \partial\Omega$.

\par\vskip 10pt
\proof
We extend the domain of $f$ to all of
$\overline{\Omega}\times \RR$ by setting
$f(x,s) = 0$ for $(x,s)\in \overline{\Omega}\times (-\infty,0)$.
We define the potential $F$ for the function $f$ as follows:
\begin{equation}
\label{e:F'=f-probl_1}
  F(x,u)\eqdef \int_0^u f(x,s) \,\mathrm{d}s =
\left\{
\begin{array}{cl}
  \int_0^u f(x,s) \,\mathrm{d}s &\quad\mbox{ if }\, 0\leq u < \infty \,;
\\
  0 &\quad\mbox{ if }\, -\infty < u < 0 \,,
\end{array}
\right.
\end{equation}
for $(x,u)\in \overline{\Omega}\times \RR$.
Hence,
\begin{math}
  f(x,s) = \frac{\partial F}{\partial u} (x,s)
\end{math}
for $(x,s)\in \overline{\Omega}\times \RR$.
Clearly, for each fixed $x\in \overline{\Omega}$,
$F(x,\,\cdot\,)\colon \RR\to \RR_+$
is a monotone increasing function, owing to
\begin{math}
  \frac{\partial F}{\partial u} (x,s) = f(x,s)\geq 0 \,.
\end{math}

Next, we obtain a nonnegative weak solution to
problem~\eqref{e:problem_1} from a global minimizer of
the energy functional
$\mathcal{E}\colon W^{1,p(x)}_0(\Omega)\to \RR$
defined by
\begin{equation}
\label{def:E(v)}
  \mathcal{E}(u)\equiv \mathcal{E}_{p(x)}(u)\eqdef
  \int_{\Omega} \frac{1}{p(x)}\, |\nabla u(x)|^{p(x)} \,\mathrm{d}x
  - \int_{\Omega} F(x,u(x)) \,\mathrm{d}x
\end{equation}
for every function $u\in W_0^{1,p(x)}(\Omega)$.
This functional is well\--defined, by the Sobolev embedding
$W_0^{1,p(x)}(\Omega) \hookrightarrow L^r(\Omega)$,
which is even compact,
and the estimate in \eqref{lim:s=infty}.
The reader is referred to the monograph by
{\sc L.\ Diening}, {\sc P.\ Harjulehto}, {\sc P.\ H\"ast\"o}, and
{\sc M.\ R\r{u}\v{z}i\v{c}ka}
\cite[{\S}8.3 and {\S}8.4]{DienHHR}
for Sobolev embeddings and their compactness.
Furthermore,
$\mathcal{E}\colon W_0^{1,p(x)}(\Omega)\to \RR$
is coercive thanks to ineq.~\eqref{lim:s=infty}
and $r\leq p_{-}$, i.e.,
\begin{equation}
\label{e:E-coerce}
  \| u\|_{ W^{1,p(x)}_0(\Omega) }\eqdef
  \|\nabla u\|_{ L^{p(x)}(\Omega) } \,\longrightarrow\, +\infty
  \quad\Longrightarrow\quad \mathcal{E}(u)\to +\infty \,.
\end{equation}
It is also weakly lower semicontinuous, by
\cite[{\S}13.2, pp.\ 412--417]{DienHHR}.
Thus, by a basic result from the calculus of variations
({\sc M. Struwe} \cite[Theorem 1.2, p.~4]{Struwe}),
$\mathcal{E}$ possesses a global minimizer
$u_0\in W_0^{1,p(x)}(\Omega)$.
Since also
$|u_0|\in W_0^{1,p(x)}(\Omega)$ with the Sobolev gradient
$\nabla |u_0| = \nabla u_0$ almost everywhere in the set
$\Omega^{+} = \{ x\in \Omega\colon u_0(x)\geq 0\}$, and
$\nabla |u_0| = - \nabla u_0$ almost everywhere in
$\Omega^{-} = \{ x\in \Omega\colon u_0(x)\leq 0\}$, we have
$| \nabla |u_0| | = |\nabla u_0|$ a.e.\ in
$\Omega = \Omega^{+}\cup \Omega^{-}$.
From this equality, combined with
$F(x,u) > 0$ for $u>0$ and $F(x,u) = 0$ for $u\leq 0$,
we deduce that
\begin{math}
  \mathcal{E}(|u_0|)\leq \mathcal{E}(u_0)
\end{math}
which shows that also $|u_0|$ is a a global minimizer for
$\mathcal{E}$ on $W_0^{1,p(x)}(\Omega)$.
This means that
\begin{math}
  \mathcal{E}(u_0)\leq \mathcal{E}(|u_0|)\leq \mathcal{E}(u_0)
\end{math}
which forces
\begin{align*}
    \mathcal{E}(u_0) =
& \int_{\Omega} \frac{1}{p(x)}\, |\nabla u_0(x)|^{p(x)} \,\mathrm{d}x
  - \int_{\Omega} F(x,u_0(x)) \,\mathrm{d}x
\\
  = \mathcal{E}(|u_0|) =
& \int_{\Omega} \frac{1}{p(x)}\,
    \bigl\vert \nabla |u_0(x)| \bigr\vert^{p(x)} \,\mathrm{d}x
  - \int_{\Omega} F(x, |u_0(x)| ) \,\mathrm{d}x \,.
\end{align*}
The arguments above yield
\begin{equation*}
    \int_{\Omega^{-}} F(x, |u_0(x)| ) \,\mathrm{d}x
  = \int_{\Omega^{-}} F(x,u_0(x)) \,\mathrm{d}x = 0 \,.
\end{equation*}
Hence, we get $u_0(x) = 0$ for a.e.\ $x\in \Omega^{-}$.
We have proved that $u_0\geq 0$ a.e.\ in $\Omega$.

We now exclude the possibility that $u_0\equiv 0$ in $\Omega$,
i.e., $u_0 = 0$ a.e.\ in $\Omega$.
Since $\mathcal{E}(0) = 0$, we only need to find a function
$u_1\in W_0^{1,p(x)}(\Omega)$ such that $\mathcal{E}(u_1) < 0$.
Then
$\mathcal{E}(u_0)\leq \mathcal{E}(u_1) < 0$ prevents the case
$u_0\equiv 0$ in $\Omega$ with $\mathcal{E}(u_0) = 0$.
To this end, choose
$\phi\in C^1_{\mathrm{c}}(\Omega)$
to be an arbitrary nonnegative $C^1$-function
with compact support in $\Omega$, $\phi\not\equiv 0$ in $\Omega$.
For $0 < t\leq 1$ we estimate
\begin{align}
\label{e:E-zero}
    \mathcal{E}(t\phi)
& = \int_{\Omega} \frac{ t^{p(x)} }{p(x)}\,
    |\nabla\phi(x)|^{p(x)} \,\mathrm{d}x
  - \int_{\Omega} F(x, t\phi(x)) \,\mathrm{d}x
\\
\nonumber
& \leq \frac{ t^{p_{-}} }{p_{-}}
  \int_{\Omega}
    |\nabla\phi(x)|^{p(x)} \,\mathrm{d}x
  - \int_{\Omega} F(x, t\phi(x)) \,\mathrm{d}x \,.
\end{align}
In order to estimate the last integral, from
ineq.~\eqref{lim:s=0} we deduce that,
given any $\eps > 0$, there is a constant $t_{\eps}\in (0,1]$
such that
\begin{equation}
\label{lim_F:s=0}
  F(x, t\phi(x))\geq \frac{1}{r\eps}\, [ t\, \phi(x)]^r
  \quad\mbox{ holds for all }\,
  (x,t)\in \overline{\Omega}\times [0,t_{\eps}] \,.
\end{equation}
We apply this estimate to ineq.~\eqref{e:E-zero}
and recall that $1 < r\leq p_{-}$, thus arriving at
\begin{align*}
    \mathcal{E}(t\phi)
& \leq \frac{t^r}{r}
  \int_{\Omega} |\nabla\phi(x)|^{p(x)} \,\mathrm{d}x
  - \frac{t^r}{r\eps}
  \int_{\Omega} \phi(x)^r \,\mathrm{d}x
\\
& = {}- \frac{t^r}{r}
    \left(
    \frac{1}{\eps} \int_{\Omega} \phi(x)^r \,\mathrm{d}x
  - \int_{\Omega} |\nabla\phi(x)|^{p(x)} \,\mathrm{d}x
    \right)
\end{align*}
for all $t\in [0,t_{\eps}]$.
Choosing $\eps > 0$ small enough, we conclude that
$\mathcal{E}(t\phi) < 0$ whenever $0 < t\leq t_{\eps}$.
In addition to $u_0\geq 0$ a.e.\ in $\Omega$,
we have proved also $u_0\not\equiv 0$ in~$\Omega$.

Since $u_0\in W_0^{1,p(x)}(\Omega)$
is a global minimizer for the functional
$\mathcal{E}\colon W_0^{1,p(x)}(\Omega)\to \RR$,
it is also a critical point for $\mathcal{E}$ and, hence,
a nonnegative weak solution to problem~\eqref{e:problem_1}
provided $u_0\in L^{\infty}(\Omega)$.

Now let $u\in W_0^{1,p(x)}(\Omega)$
be any nonnegative critical point for $\mathcal{E}$,
$u\not\equiv 0$ in~$\Omega$.
This means that $u$ is a weak solution to problem~\eqref{e:problem_1}
in the sense of
{\sc X.\ Fan} and {\sc D.\ Zhao} \cite[Def.\ 4.1, p.~311]{Fan-Zhao}.
We may apply their regularity result
\cite[Theorem 4.1, p.~312]{Fan-Zhao}
(and its proof)
to conclude that $u\in L^{\infty}(\Omega)$.
This means that $u$ is a nonnegative weak solution
to problem~\eqref{e:problem_1}
also in the sense of our Definition~\ref{def-weak_sol} above.
Moreover, we get $u\in C^{0,\beta'}(\Omega)$
for some $\beta'\in (0, \alpha)$, by
\cite[Theorem 4.2, p.~315]{Fan-Zhao}.
Furthermore, thanks to our Hypothesis {\rm ($\mathbf{p}$)} on~$p$, i.e.,
$p\in C^{0, \alpha_1}(\overline{\Omega})$ for some $\alpha_1\in (0,1)$,
we may apply a stronger regularity result due to
{\sc X.-L.\ Fan} \cite[Theorem 1.2, p.~400]{Fan-JDE} to obtain
$u\in C^{1,\beta}(\overline{\Omega})$ for some $\beta\in (0, \alpha)$.
Finally, we apply the strong maximum principle and
the Hopf boundary point lemma, respectively, from
{\sc Q.\ Zhang} \cite[Theorems 1.1 and 1.2, p.~26]{Zhang}
to conclude that both inequalities claimed in \eqref{Hopf:problem_1}
are valid.

Clearly, the global minimizer $u_0\in W_0^{1,p(x)}(\Omega)$
for the functional $\mathcal{E}$ obtained above
enjoys analogous regularity and positivity properties as does $u$.
As a simple consequence, both ratios
$u / u_0$ and $u_0 / u$ are continuous positive functions
over the domain $\Omega$ and can be extended to
positive continuous functions over the closure $\overline{\Omega}$,
by l'Hospital's rule,
\begin{equation}
\label{e:l'Hospital}
    \lim_{x\to x_0} \frac{u(x)}{u_0(x)}
  = \lim_{t\to 0+}
    \frac{ u  \left( x_0 - t \boldsymbol{\nu} (x_0) \right) }%
         { u_0\left( x_0 - t \boldsymbol{\nu} (x_0) \right) }
  = \frac{\partial u  }{\partial\boldsymbol{\nu}} (x_0) \bigg\slash
    \frac{\partial u_0}{\partial\boldsymbol{\nu}} (x_0)
  > 0 \,,
\end{equation}
where $x_0\in \partial\Omega$ is an arbitrary boundary point and
$x\in \Omega$ ranges inside $\Omega$ near $x_0$, e.g.,
$x = x_0 - t \boldsymbol{\nu} (x_0)$ with $t > 0$ small enough.
The last ratio,
\begin{math}
  x_0\mapsto
  \frac{\partial u  }{\partial\boldsymbol{\nu}} (x_0) \Big\slash
  \frac{\partial u_0}{\partial\boldsymbol{\nu}} (x_0) \,,
\end{math}
being positive and continuous over the compact boundary
$\partial\Omega$, we conclude that both ratios,
$u / u_0$ and $u_0 / u$, can be extended to
positive continuous functions over the closure $\overline{\Omega}$.
Consequently, both ratios are bounded.
We apply our Theorem~\ref{thm-Diaz-Saa}
(the {\sc D{\'\i}az} and {\sc Saa} inequality)
to arrive at the uniqueness of a nonnegative and nontrivial weak solution
$u\in W^{1,p(x)}_0(\Omega)\cap L^{\infty}(\Omega)$
to problem~\eqref{e:problem_1}, i.e., $u = u_0$, as follows.
Setting $w_1 = u$ and $w_2 = u_0$ in Theorem~\ref{thm-Diaz-Saa},
the left\--hand side of ineq.~\eqref{in_A:Diaz-Saa} becomes
\begin{align}
\nonumber
& \int_{\Omega}
  \left(
    \frac{ {}- \Delta_{p(x)} u   }{ u  (x)^{r-1} }
  - \frac{ {}- \Delta_{p(x)} u_0 }{ u_0(x)^{r-1} }
  \right)
  (u(x)^r - u_0(x)^r) \,\mathrm{d}x
\\
\label{in_A:Diaz-Saa:u,u_0}
  =
& \int_{\Omega}
  \left(
    \frac{ f(x, u  (x)) }{ u  (x)^{r-1} }
  - \frac{ f(x, u_0(x)) }{ u_0(x)^{r-1} }
  \right)
  (u(x)^r - u_0(x)^r) \,\mathrm{d}x\leq 0 \,,
\end{align}
since the function
\begin{math}
  s\mapsto f(x,s) / s^{r-1} \colon (0,\infty)\to \RR_+
\end{math}
is strictly monotone decreasing for every $x\in \Omega$,
by Hypothesis~{\bf (f2)}.
However, by ineq.~\eqref{in_A:Diaz-Saa},
precisely the opposite inequality ``$\geq$'' must be valid.
We conclude that the equality in \eqref{in_A:Diaz-Saa:u,u_0} above
must hold.
That is possible only if
$u(x) = u_0(x)$ at almost every point $x\in \Omega$,
by Hypothesis~{\bf (f2)}, i.e.,
$u\equiv u_0$ in $\Omega$, by the regularity derived above.

Our proof of Theorem~\ref{thm-problem_1} is now complete.
\qed
\par\vskip 10pt

\begin{remark}\label{rem-problem_1}\nopagebreak
\begingroup\rm
In our proof of Theorem~\ref{thm-problem_1} above we have proved that
any nonnegative critical point
$u\in W_0^{1,p(x)}(\Omega)$ for the energy functional
$\mathcal{E}\colon W^{1,p(x)}_0(\Omega)\to \RR$
defined by eq.~\eqref{def:E(v)}
must be bounded, i.e.,
$u\in W_0^{1,p(x)}(\Omega)\cap L^{\infty}(\Omega)$.
Hence, $u$ is a nonnegative weak solution
to problem~\eqref{e:problem_1}
in the sense of our Definition~\ref{def-weak_sol} above.
The decisive argument here is the regularity result in
{\sc X.\ Fan} and {\sc D.\ Zhao} \cite[Theorem 4.1, p.~312]{Fan-Zhao}.
\endgroup
\end{remark}
\par\vskip 10pt

\begin{example}\label{exam-appl_1}\nopagebreak
\begingroup\rm
{\rm (a)}$\;$
A typical example of the function $f$ satisfying all conditions
in Theorem~\ref{thm-problem_1} is
$f(x,s) = h(x)\, s^{q(x)-1}$, with a positive function
$h\in C(\overline{\Omega})$ and
$q\in C(\overline{\Omega})$ such that
\begin{equation*}
  1\leq q(x)\leq q_{+}\eqdef \sup_{\Omega} q(x) < r = p_-
    \quad\mbox{ for every }\, x\in \overline{\Omega} \,.
\end{equation*}
Consequently, $f$ satisfies Hypotheses {\bf (f1)} -- {\bf (f3)}
with $r = p_-$.

{\rm (b)}$\;$
Our condition on $r$, i.e.,
$q(x) < r = p_-$ for every $x\in \overline{\Omega}$,
is trivially sharp in the following sense:
If $q(x)\equiv r\equiv p(x)$ is a constant in $\overline{\Omega}$,
$1 < r < \infty$, and $h(x)\equiv \lambda_{1,r}(\Omega)$
is the first eigenvalue of the ``positive'' $r$-Laplacian
\begin{math}
  -\Delta_r u = {}- \mathrm{div} (|\nabla u|^{r-2} \nabla u)
\end{math}
in $\Omega$ with zero Dirichlet boundary conditions, and
$\varphi_{1,r}\in W_0^{1,r}(\Omega)$
denotes the associated first eigenfunction normalized by
$\varphi_{1,r}(x) > 0$ in $\Omega$ and
$\int_{\Omega} \varphi_{1,r}(x)^r \,\mathrm{d}x = 1$,
then any nonnegative multiple $t\varphi_{1,r}$ ($t\in \RR_+$)
is a nonnegative weak solution to problem~\eqref{e:problem_1};
hence, this problem admits an infinite number of solutions.

In contrast, if
$0\leq h(x) < \lambda_{1,r}(\Omega)$ holds for all $x\in \Omega$,
then the variational characterization of the first eigenvalue
$\lambda_{1,r}(\Omega)$ by the Rayleigh quotient
\begin{equation*}
  \lambda_{1,r}(\Omega) = \min_{u\not\equiv 0}
  \frac{ \int_{\Omega} |\nabla u(x)|^r \,\mathrm{d}x }%
       { \int_{\Omega} |u(x)|^r        \,\mathrm{d}x }
\end{equation*}
with every minimizer $\varphi\in W_0^{1,r}(\Omega)$
taking the form $\varphi = t\varphi_{1,r}$ for some
$t\in \RR\setminus \{ 0\}$,
leaves only the trivial zero solution to problem~\eqref{e:problem_1}
with $f(x,s) = h(x)\, s^{r-1}$.

{\rm (c)}$\;$
Furthermore, if
$1\leq q(x)\leq r = p_{-} = \inf_{x\in \Omega} p(x)$
holds for every $x\in \Omega$, with
$p(x) > r$ for all $x\in \Omega_0$ in a subset
$\Omega_0\subset \Omega$ of positive Lebesgue measure,
then problem~\eqref{e:problem_1} possesses
at most one nonnegative and nontrivial weak solution
$u\in W^{1,p(x)}_0(\Omega)\cap L^\infty(\Omega)$.
Any such weak solution $u$ belongs to the class
$C^{1,\beta}(\overline{\Omega})$, for some $\beta\in (0, \alpha)$,
and satisfies also the Hopf maximum principle \eqref{Hopf:problem_1}.
Of course, $u = 0$ on the boundary $\partial\Omega$.
Hence, $u$ is also a positive weak solution.
This claim follows easily from the fact that the reaction function
$f(x,s) = h(x)\, s^{q(x)-1}$ has the following properties:
\begin{itemize}
\item[{\bf (f1')}]
$\;$
$f\colon \overline{\Omega}\times \RR_+\to \RR_+$
is a nonnegative continuous function satisfying
$f(x,0) = 0$ for all $x\in \Omega$ with $q(x) > 1$.
\item[{\bf (f2')}]
$\;$
The function
\begin{equation*}
  s\,\longmapsto\, \frac{ f(x,s) }{ s^{r-1} }
  = \frac{h(x)}{s^{r - q(x)}}
  \,\colon (0,\infty)\to \RR_+
\end{equation*}
is strictly monotone decreasing for every
$x\in \Omega_1 = \{ x\in \Omega\colon q(x) < r\}$,
while being ${}= h(x)$ for all
\begin{math}
  x\in \Omega\setminus \Omega_1 = \{ x\in \Omega\colon q(x) = r\} .
\end{math}
\end{itemize}
Recall that $h>0$ in all of $\Omega$.
Consequently, if there were two distinct positive weak solutions, say,
$u_0$ and $u$ as in our proof of Theorem~\ref{thm-problem_1} above,
then inequality \eqref{in_A:Diaz-Saa:u,u_0} would force
$u(x) = u_0(x)$ for every $x\in \Omega_1$.
Moreover, by Theorem~\ref{thm-Diaz-Saa}, Part~{\rm (b)},
we have even $u(x) = u_0(x)$ for every $x\in \Omega$,
thanks to $p(x) > r$ for all $x\in \Omega_0$.

{\rm (d)}$\;$
Similarly to case {\rm (c)} above, if
$1\leq q(x)\leq r = p_{-}\equiv p(x)$
holds for every $x\in \Omega$, with
$p(x)\equiv r\in (1,\infty)$ being a constant and
$q(x) < r$ for all $x\in \Omega_1$ in a subset
$\Omega_1\subset \Omega$ of positive Lebesgue measure,
then problem~\eqref{e:problem_1} possesses
at most one nonnegative and nontrivial weak solution
$u\in W^{1,p(x)}_0(\Omega)\cap L^\infty(\Omega)$.
The reasoning for this is similar as in case {\rm (c)}:
First, we may take
$\Omega_1 = \{ x\in \Omega\colon q(x) < r\}$;
its Lebesgue measure is ${}> 0$.
Again, inequality \eqref{in_A:Diaz-Saa:u,u_0} forces
$u(x) = u_0(x)$ for every $x\in \Omega_1$.
Moreover, by Theorem~\ref{thm-Diaz-Saa}, Part~{\rm (a)},
we have even $u(x) = c\cdot u_0(x)$ for almost every $x\in \Omega$,
where $c\in (0,\infty)$ is a constant.
But $\Omega_1\subset \Omega$ has positive Lebesgue measure
which yields $c=1$.
The uniqueness result follows.

{\rm (e)}$\;$
Finally, in
{\sc M.\ Mih\u{a}ilescu} and {\sc V. R\u{a}dulescu} \cite{Mihai-Radu},
the nonuniqueness of weak solutions is established in case
\begin{math}
  1 < \min_{ x\in \overline{\Omega} } q(x) < p_-
    < \max_{ x\in \overline{\Omega} } q(x)
\end{math}
(see also some other related results in case $p(x) = q(x)$
 in {\sc X. Fan}, {\sc Q. Zhang}, and {\sc D. Zhao}
 \cite{Fan-Zhang-Zhao}).
\endgroup
\end{example}
\par\vskip 10pt

\begin{remark}\label{rem-problem_1-rem}\nopagebreak
\begingroup\rm
{\rm (i)}$\;$
Theorem~\ref{thm-problem_1} solves the open problem raised in
{\sc X. Fan} \cite{Fan-Math-N} (see Remark 2.3 on p.~1443)
and improves the uniqueness results given for
problem~\eqref{e:problem_1} in \cite{Fan-Math-N}.

{\rm (ii)}$\;$
The uniqueness property does not hold for solutions with changing sign,
even if $p$ is a constant.
For more details, we refer to examples exhibiting
two distinct critical points of the energy functional
${\mathcal E}_\lambda$, with $\lambda > 0$, defined on
$W^{1,p}_0(\Omega)$ by
\begin{eqnarray*}
  {\mathcal E}_\lambda(u)\eqdef
    \frac{1}{p}\int_\Omega \vert\nabla u\vert^p \,\mathrm{d}x
  - \frac{\lambda}{p}\int_\Omega \vert u\vert^p \,\mathrm{d}x
  - \int_{\Omega} f(x)\, u \,\mathrm{d}x \,.
\end{eqnarray*}
Such examples were constructed in
{\sc M.~A.\ del Pino}, {\sc M.\ Elgueta}, and {\sc R.~F.\ Man\'asevich}
\cite[Eq.\ (5.26) on p.~12]{PiElMa} for $2 < p < \infty$ and in
{\sc J. Fleckinger}, {\sc J. Hern\'andez}, and {\sc P. Tak\'a\v{c}}
\cite[Example~2 on p.~148]{FleckHTT} for $1 < p < 2$.

{\rm (iii)}$\;$
The easiest problem of type \eqref{e:problem_1},
with the right\--hand side $f(x,u)\equiv f(x)$
being independent from the unknown function $u = u(x)$,
$f\in L^{\infty}(\Omega)$,
can be treated in a similar way as in Theorem~\ref{thm-problem_1};
one has to take $r=1$ in the proof, particularly in
Theorem~\ref{thm-Diaz-Saa} when applying it to an anlogue of
ineq.~\eqref{in_A:Diaz-Saa:u,u_0}.
Then this inequality is actually an equality with
the right\--hand side $= 0$.
The uniqueness of a weak solution to problem
\begin{equation*}
  {}- \Delta_{p(x)} u = f(x) \;\mbox{ in }\, \Omega \,;\quad
  u = 0 \;\mbox{ on }\, \partial\Omega \,,\quad
  u > 0 \mbox{ in }\, \Omega \,,
\end{equation*}
then follows by Theorem~\ref{thm-Diaz-Saa}, Part~{\rm (a)}.
But this uniqueness result is valid for any
(possibly sign\--changing) weak solution $u\in W^{1,p(x)}_0(\Omega)$,
by a classical argument that takes advantage of
the strict convexity of the functional
\begin{math}
  \mathcal{W}_A\equiv \mathcal{W}_{A,p(x),1}
\end{math}
on $W^{1,p(x)}_0(\Omega)$ ($r=1$); see, e.g.,
\cite[Theorem 13.3.3, p.~418]{DienHHR} or
\cite[Theorem 4.3, p.~1848]{Fan-Zhang-1}.
\endgroup
\end{remark}
\par\vskip 10pt

Our second example is the following simple generalization of
problem~\eqref{e:problem_1}:
\begin{equation}
\label{e:problem_2}
\left\{
\begin{alignedat}{2}
  {}- \Delta_{p(x)} u + g(x,u)
& {}= f(x,u)
&&  \quad\mbox{ in }\, \Omega \,;
\\
  u & {}= 0
&&  \quad\mbox{ on }\, \partial\Omega \,,
    \quad u > 0 \mbox{ in }\, \Omega \,.
\end{alignedat}
\right.
\end{equation}
Here, we have a new monotone nonlinear operator on the left\--hand side,
${}- \Delta_{p(x)} u + g(x,u)$, whose homogeneity properties
with respect to the function $u$ are similar to those of
\hfil\break
\begin{math}
  {}- \Delta_{p(x)} u =
  {}- \mathrm{div} \left( |\nabla u|^{p(x)-2} \nabla u\right) .
\end{math}
We recall that $p\colon \overline{\Omega}\to (1,\infty)$
is a continuous function, such that
it satisfies Hypothesis ($\mathbf{p}$)
together with inequalities \eqref{e:p-<p<p+}, where
$r\in \RR$ is a given constant, $1 < r\leq p_{-}$.
The function
$f\colon \overline{\Omega}\times \RR_+\to \RR_+$
is assumed to satisfy all Hypotheses {\bf (f1)} -- {\bf (f3)}.

\par\vskip 10pt
%
We impose the following hypotheses on the function $g$:
\begin{itemize}
\item[{\bf (g1)}]
$\;$
$g\colon \overline{\Omega}\times \RR_+\to \RR_+$
is a nonnegative continuous function such that
$g(x,0) = 0$ for all $x\in \Omega$ and
$g(x,s) > 0$ for all $(x,s)\in \Omega\times (0,\infty)$.
\item[{\bf (g2)}]
$\;$
The function
\begin{math}
  s\,\longmapsto\, g(x,s) / s^{r-1} \colon (0,\infty)\to \RR_+
\end{math}
is monotone increasing for every $x\in \Omega$,
but {\em not\/} necessarily {\em strictly\/} monotone increasing.
\item[{\bf (g3)}]
$\;$
The following limit is uniform with respect to $x\in \Omega$:
\begin{equation*}
    \limsup_{s\to +\infty}
  \frac{ g(x,s) }{ s^{m(x)-1} } \leq C\equiv \mathrm{const} < \infty
  \quad\mbox{ for all }\, x\in \Omega \,,
\end{equation*}
where $m\colon \overline{\Omega}\to \RR_+$ is
some suitable continuous function that satisfies
$1 < m(x) < p^{\ast}(x)$, where
\begin{equation*}
  p^{\ast}(x)\eqdef
  \left\{
\begin{alignedat}{2}
& \genfrac{}{}{}1{N p(x)}{N - p(x)} \quad &&\mbox{ if }\; p(x) < N \,;
\\
& +\infty \quad &&\mbox{ if }\; p(x)\geq N \,.
\end{alignedat}
\right.
\end{equation*}
\end{itemize}
%
\par\vskip 10pt

The authors in \cite[{\S}8.3, pp.\ 265--272]{DienHHR}
call $p^{\ast}(x)\in [1, +\infty]$
the {\em Sobolev conjugate exponent\/} and prove the Sobolev embedding
\begin{math}
  W_0^{1,p(\cdot)}(\Omega) \hookrightarrow L^{p^{\ast}(\cdot)}(\Omega)
\end{math}
for
\begin{math}
  p_{+} = \sup_{\Omega} p(x)\break
  < N
\end{math}
(\cite[{\S}8.3, Theorem~8.3.1, p.~265]{DienHHR})
under the additional regularity hypothesis on $p(x)$ requiring
$p\in \mathcal{P}^{\log}(\Omega)$, cf.\
\cite[{\S}4.1, Def.~4.1.4, p.~101]{DienHHR}, i.e.,
$1/p(x)$ is globally log\--H\"older\--continuous in $\Omega$.
This additional hypothesis (log\--H\"older continuity)
is always satisfied in our situation, provided
$p\colon \overline{\Omega}\to (1,\infty)$
is a continuous function that obeys our hypotheses above, i.e.,
$p$ satisfies Hypothesis ($\mathbf{p}$)
together with inequalities \eqref{e:p-<p<p+},
where we now assume also $p_{+} < N$, in addition to $1 < r\leq p_{-}$.
%
\par\vskip 10pt

It is worth of noticing that
Hypotheses {\bf (g1)} and~{\bf (g2)} imply
\begin{itemize}
\item[{\bf (g2')}]
$\;$
Also
\begin{math}
  s\,\longmapsto\, g(x,s)\colon \RR_+\to \RR_+
\end{math}
is a {\em strictly\/} monotone increasing function for every 
$x\in \Omega$.
\end{itemize}
Moreover, Hypotheses {\bf (g1)} and~{\bf (g2)} combined entail
\begin{equation}
\label{g-lim:s=0}
  g(x,s)\leq C_0\, s^{r-1}
  \quad\mbox{ holds for all }\,
  (x,s)\in \overline{\Omega}\times [0,s_0] \,.
\end{equation}
Here, $s_0\in [1,\infty)$ is an arbitrary number and
\begin{equation*}
  C_0 = C_0(s_0) = \frac{ \sup_{x\in \Omega} g(x,s_0) }{s_0^{r-1}}
  < \infty
\end{equation*}
is a positive constant depending solely on $s_0$.

We remark that Hypotheses {\bf (g1)} and~{\bf (g3)}
limit the asymptotic behavior of $g(x,\,\cdot\,)$
as $s\to +\infty$ as follows:
Given any $\eps > 0$, there is a constant $C_{\eps}'\in (0,\infty)$
such that
\begin{equation}
\label{lim:s=infty2}
  0\leq g(x,s)\leq (C+\eps)\, s^{m(x)-1} + C_{\eps}'
  \quad\mbox{ holds for all }\,
  (x,s)\in \overline{\Omega}\times \RR_+ \,.
\end{equation}

In analogy with our Definition~\ref{def-weak_sol}
adapted to problem \eqref{e:problem_1},
we define the notion of a {\it nonnegative\/} weak solution to
problem \eqref{e:problem_2} as follows:

\begin{definition}\label{def-weak_sol_2}\nopagebreak
\begingroup\rm
A nonnegative function
$u\in W^{1,p(x)}_0(\Omega)\cap L^\infty(\Omega)$
is called a {\em\bfseries nonnegative weak solution\/}
of problem \eqref{e:problem_2} if,
for every test function $\phi\in W^{1,p(x)}_0(\Omega)$,
the following equation holds,
\begin{equation}
\label{def:problem_2}
\begin{aligned}
& \int_{\Omega}
    \vert\nabla u(x)\vert^{p(x)-2}\nabla u(x)\cdot \nabla\phi(x)
   \,\mathrm{d}x 
  + \int_{\Omega} g(x,u(x))\, \phi(x) \,\mathrm{d}x
\\
& {}
  = \int_{\Omega} f(x,u(x))\, \phi(x) \,\mathrm{d}x \,.
\end{aligned}
\end{equation}

If\/ $u$ satisfies also $u > 0$ throughout $\Omega$, we call $u$
a {\em\bfseries positive weak solution\/}.
\endgroup
\end{definition}
\par\vskip 10pt

Problem~\eqref{e:problem_2}
fits into a more general class of variational problems treated in
{\sc X.\ Fan} and {\sc D.\ Zhao} \cite[Eq.\ 4.1, p.~310]{Fan-Zhao}.
However, the authors are interested only in some standard
regularity properties of weak solutions, like
(local and global) boundedness and H\"older continuity
(\cite[Sect.~4, pp.\ 310--317]{Fan-Zhao}).

We now generalize the existence and uniqueness result in
Theorem~\ref{thm-problem_1}
to the boundary value problem~\eqref{e:problem_2}
for positive weak solutions.

\begin{theorem}\label{thm-problem_2}
Under the {\rm Hypotheses\/}
{\rm ($\boldsymbol{\Omega}$)}, {\rm ($\mathbf{p}$)},
{\bf (f1)} -- {\bf (f3)}, and\/
{\bf (g1)} -- {\bf (g3)},
problem~\eqref{e:problem_2} possesses
a unique nonnegative and nontrivial weak solution
$u\in W^{1,p(x)}_0(\Omega)\cap L^\infty(\Omega)$.
This solution belongs to the class
$C^{1,\beta}(\overline{\Omega})$, for some $\beta\in (0, \alpha)$,
and satisfies also the Hopf maximum principle
\eqref{Hopf:problem_1},
\begin{equation*}
  u(x) > 0 \;\mbox{ for all }\, x\in \Omega
  \quad\mbox{ and }\quad
  \frac{\partial u}{\partial\boldsymbol{\nu}} (x) < 0
    \;\mbox{ for all }\, x\in \partial\Omega \,.
\end{equation*}
Of course, $u = 0$ on the boundary $\partial\Omega$.
Hence, $u$ is also a positive weak solution.
\end{theorem}

\par\vskip 10pt
\proof
First, let us recall that the potential $F$ for the function $f$
has been defined in eq.~\eqref{e:F'=f-probl_1}.
We define the potential $G$ for the function $g$
in a similar way:
First, we extend the domain of $g$ to all of
$\overline{\Omega}\times \RR$ by setting
$g(x,s) = 0$ for $(x,s)\in \overline{\Omega}\times (-\infty,0)$.
Then we define the potential $G$ for the function $g$ by
\begin{equation}
\label{e:G'=g-probl_2}
  G(x,u)\eqdef \int_0^u g(x,s) \,\mathrm{d}s =
\left\{
\begin{array}{cl}
  \int_0^u g(x,s) \,\mathrm{d}s &\quad\mbox{ if }\, 0\leq u < \infty \,;
\\
  0 &\quad\mbox{ if }\, -\infty < u < 0 \,,
\end{array}
\right.
\end{equation}
for $(x,u)\in \overline{\Omega}\times \RR$.
Hence,
\begin{math}
  g(x,s) = \frac{\partial G}{\partial u} (x,s)
\end{math}
for $(x,s)\in \overline{\Omega}\times \RR$.
Clearly, for each fixed $x\in \overline{\Omega}$,
$G(x,\,\cdot\,)\colon \RR\to \RR_+$
is a monotone increasing function, owing to
\begin{math}
  \frac{\partial G}{\partial u} (x,s) = g(x,s)\geq 0 \,.
\end{math}

Again, we obtain a nonnegative weak solution to
problem~\eqref{e:problem_2} from a global minimizer of
the energy functional
$\hat{\mathcal{E}}\colon W^{1,p(x)}_0(\Omega)\to \RR$
defined by
\begin{equation}
\label{def:E_2(v)}
\begin{aligned}
  \hat{\mathcal{E}}(u)\equiv
  \hat{\mathcal{E}}_{p(x)}(u)\eqdef
& \int_{\Omega} \frac{1}{p(x)}\, |\nabla u(x)|^{p(x)} \,\mathrm{d}x
\\
& {}
  + \int_{\Omega} G(x,u(x)) \,\mathrm{d}x
  - \int_{\Omega} F(x,u(x)) \,\mathrm{d}x
\end{aligned}
\end{equation}
for every function $u\in W_0^{1,p(x)}(\Omega)$.
By our proof of Theorem~\ref{thm-problem_1},
the first and last summands on the right\--hand side of
eq.~\eqref{def:E_2(v)} are well\--defined.
The same is true of the second summand, thanks to inequalities
\eqref{g-lim:s=0} and \eqref{lim:s=infty2}
supplemented by the Sobolev embedding
\begin{math}
  W_0^{1,p(\cdot)}(\Omega) \hookrightarrow L^{p^{\ast}(\cdot)}(\Omega)
\end{math}
for $p_{+} < N$
(\cite[{\S}8.3, Theorem~8.3.1, p.~265]{DienHHR}).

By the standard properties of the ``smaller'' functional
\begin{math}
  \mathcal{E}(u) =
  \hat{\mathcal{E}}(u) {}- \hfil\break
  \int_{\Omega} G(x,u(x)) \,\mathrm{d}x
\end{math}
defined in eq.~\eqref{def:E(v)}, that have been verified
in the proof of Theorem~\ref{thm-problem_1} above,
also our present functional
$\hat{\mathcal{E}}\colon W_0^{1,p(x)}(\Omega)\to \RR$
is coercive thanks to ineq.~\eqref{lim:s=infty}
and $r\leq p_{-}$, i.e.,
it satisfies an analogue of \eqref{e:E-coerce},
\begin{equation*}
  \| u\|_{ W^{1,p(x)}_0(\Omega) }\eqdef
  \|\nabla u\|_{ L^{p(x)}(\Omega) } \,\longrightarrow\, +\infty
  \quad\Longrightarrow\quad \hat{\mathcal{E}}(u)\to +\infty \,.
\end{equation*}
It is also weakly lower semicontinuous, by
\cite[{\S}13.2, pp.\ 412--417]{DienHHR}.
Consequently, a basic result from the calculus of variations yields
the existence of a global minimizer
$u_0\in W_0^{1,p(x)}(\Omega)$ for $\mathcal{E}$.
We claim that $u_0\geq 0$ a.e.\ in $\Omega$.
Clearly, also its positive part,
$u_0^{+}\eqdef \max\{ u_0, 0\} \geq 0$,
is in $W_0^{1,p(x)}(\Omega)$ and thus satisfies
$\hat{\mathcal{E}}(u_0^{+}) \geq \hat{\mathcal{E}}(u_0)$.
Denoting
$\Omega^{+} = \{ x\in \Omega\colon u_0(x)\geq 0\}$ and
$\Omega^{-} = \{ x\in \Omega\colon u_0(x)\leq 0\}$, we calculate
\begin{align*}
    \hat{\mathcal{E}}(u_0)
& {}
  = \int_{\Omega} \frac{1}{p(x)}\, |\nabla u_0(x)|^{p(x)} \,\mathrm{d}x
\\
& {}
  + \int_{\Omega} G(x,u_0(x)) \,\mathrm{d}x
  - \int_{\Omega} F(x,u_0(x)) \,\mathrm{d}x \,,
\end{align*}
\begin{align*}
    \hat{\mathcal{E}}(u_0)
& {}
  = \int_{ \Omega^{+} }
    \frac{1}{p(x)}\, |\nabla u_0(x)|^{p(x)} \,\mathrm{d}x
  + \int_{ \Omega^{-} }
    \frac{1}{p(x)}\, |\nabla u_0(x)|^{p(x)} \,\mathrm{d}x
\\
& {}
  + \int_{ \Omega^{+} } G(x,u_0(x)) \,\mathrm{d}x
  - \int_{ \Omega^{+} } F(x,u_0(x)) \,\mathrm{d}x
\\
& {}
  = \hat{\mathcal{E}}(u_0^{+})
  + \int_{ \Omega^{-} }
    \frac{1}{p(x)}\, |\nabla u_0(x)|^{p(x)} \,\mathrm{d}x
\\
& {}
  \geq \hat{\mathcal{E}}(u_0)
  + \int_{ \Omega^{-} }
    \frac{1}{p(x)}\, |\nabla u_0(x)|^{p(x)} \,\mathrm{d}x
  \geq \hat{\mathcal{E}}(u_0) \,.
\end{align*}
These inequalities force
$\nabla u_0^{-}(x) = {}- \nabla u_0(x)\equiv 0$
a.e.\ in $\Omega^{-}$, whence
$u_0^{-}(x)\equiv 0$ a.e.\ in~$\Omega^{-}$.
We have proved $u_0\geq 0$ a.e.\ in $\Omega$ as claimed.

In order to exclude the possibility that $u_0\equiv 0$ in $\Omega$,
we now construct a function
$u_1\in W_0^{1,p(x)}(\Omega)$ such that
$\hat{\mathcal{E}}(u_1) < 0 = \hat{\mathcal{E}}(0)$.
First, we take an arbitrary nonnegative $C^1$-function
$\phi\in C^1_{\mathrm{c}}(\Omega)$
with compact support in $\Omega$, $\phi\not\equiv 0$ in $\Omega$.
For $0 < t\leq 1$ we estimate
\begin{align}
\label{e_hat:E-zero}
    \hat{\mathcal{E}}(t\phi)
& = \int_{\Omega} \frac{ t^{p(x)} }{p(x)}\,
    |\nabla\phi(x)|^{p(x)} \,\mathrm{d}x
  + \int_{\Omega} G(x, t\phi(x)) \,\mathrm{d}x
  - \int_{\Omega} F(x, t\phi(x)) \,\mathrm{d}x
\\
\nonumber
& \leq \frac{ t^{p_{-}} }{p_{-}}
  \int_{\Omega}
    |\nabla\phi(x)|^{p(x)} \,\mathrm{d}x
  + \int_{\Omega} G(x, t\phi(x)) \,\mathrm{d}x
  - \int_{\Omega} F(x, t\phi(x)) \,\mathrm{d}x \,.
\end{align}
We estimate the difference of the last two integrals as follows.
We combine inequalities \eqref{lim:s=0} and \eqref{g-lim:s=0}
to deduce that,
given any $\eps > 0$ small enough, $\eps < 1/C_0$,
there is a constant $t_{\eps}'\in (0,1]$ such that
\begin{equation}
\label{lim_G:s=0}
\begin{aligned}
& F(x, t\phi(x)) - G(x, t\phi(x))
  \geq \frac{1}{r\eps}\, [ t\, \phi(x)]^r
     - \frac{C_0}{r}\, [ t\, \phi(x)]^r
\\
& {}
  = \frac{\eps^{-1} - C_0}{r}\, [ t\, \phi(x)]^r
  \quad\mbox{ holds for all }\,
  (x,t)\in \overline{\Omega}\times [0,t_{\eps}'] \,.
\end{aligned}
\end{equation}
Applying this estimate to ineq.~\eqref{e_hat:E-zero}
we arrive at
\begin{align*}
    \hat{\mathcal{E}}(t\phi)
& \leq \frac{t^r}{r}
  \int_{\Omega} |\nabla\phi(x)|^{p(x)} \,\mathrm{d}x
  - \frac{t^r}{r}\, (\eps^{-1} - C_0)
  \int_{\Omega} \phi(x)^r \,\mathrm{d}x
\\
\nonumber
& {}
  = {}- \frac{t^r}{r}
    \left(
    (\eps^{-1} - C_0) \int_{\Omega} \phi(x)^r \,\mathrm{d}x
    - \int_{\Omega} |\nabla\phi(x)|^{p(x)} \,\mathrm{d}x
    \right)
\end{align*}
for all $t\in [0,t_{\eps}']$.
Choosing $\eps > 0$ small enough, we conclude that
$\hat{\mathcal{E}}(t\phi) < 0$ whenever $0 < t\leq t_{\eps}'$.
In addition to $u_0\geq 0$ a.e.\ in $\Omega$,
we have proved also $u_0\not\equiv 0$ in~$\Omega$, thanks to
\begin{math}
  \hat{\mathcal{E}}(u_0) \leq
  \hat{\mathcal{E}}(t\phi) < 0 = \hat{\mathcal{E}}(0) \,.
\end{math}

Since $u_0\in W_0^{1,p(x)}(\Omega)$
is a global minimizer for the functional
$\hat{\mathcal{E}}\colon W_0^{1,p(x)}(\Omega)\to \RR$,
it is also a critical point for $\hat{\mathcal{E}}$ and, hence,
a nonnegative weak solution to problem~\eqref{e:problem_2}
provided $u_0\in L^{\infty}(\Omega)$.

Now let $u\in W_0^{1,p(x)}(\Omega)$
be any nonnegative critical point for $\hat{\mathcal{E}}$,
$u\not\equiv 0$ in~$\Omega$.
This means that $u$ is a weak solution to problem~\eqref{e:problem_2}
in the sense of
{\sc X.\ Fan} and {\sc D.\ Zhao} \cite[Def.\ 4.1, p.~311]{Fan-Zhao}.
We may apply their regularity result
\cite[Theorem 4.1, p.~312]{Fan-Zhao}
(and its proof)
to conclude that $u\in L^{\infty}(\Omega)$.
This means that $u$ is a nonnegative weak solution
to problem~\eqref{e:problem_1}
also in the sense of our Definition~\ref{def-weak_sol} above.
By another result in
\cite[Theorem 4.4, p.~317]{Fan-Zhao},
$u$ is even H\"older\--continuous in $\overline{\Omega}$,
$u\in C^{0,\beta'}(\overline{\Omega})$
for some $\beta'\in (0, \alpha)$.

The regularity property
$u\in C^{1,\beta}(\overline{\Omega})$ for some $\beta\in (0, \alpha)$
and l'Hospital's rule \eqref{e:l'Hospital}
are obtained by the same arguments as in the proof of
Theorem~\ref{thm-problem_1} above.
In particular, the continuity and boundedness of both ratios,
$u / u_0$ and $u_0 / u$, in the closure $\overline{\Omega}$ follows.
Thus, it remains to apply our Theorem~\ref{thm-Diaz-Saa}
(the {\sc D{\'\i}az} and {\sc Saa} inequality)
to arrive at the uniqueness of a nonnegative and nontrivial weak solution
$u\in W^{1,p(x)}_0(\Omega)\cap L^{\infty}(\Omega)$
to problem~\eqref{e:problem_2}, i.e., $u = u_0$.

Setting $w_1 = u$ and $w_2 = u_0$ in Theorem~\ref{thm-Diaz-Saa},
the left\--hand side of ineq.~\eqref{in_A:Diaz-Saa} becomes
\begin{align*}
& \int_{\Omega}
  \left(
    \frac{ {}- \Delta_{p(x)} u   }{ u  (x)^{r-1} }
  - \frac{ {}- \Delta_{p(x)} u_0 }{ u_0(x)^{r-1} }
  \right)
  (u(x)^r - u_0(x)^r) \,\mathrm{d}x
\\
&
\begin{aligned}
= \int_{\Omega}
& \left[
  \left(
    \frac{ f(x, u  (x)) }{ u  (x)^{r-1} }
  - \frac{ f(x, u_0(x)) }{ u_0(x)^{r-1} }
  \right)
  \right.
\\
& \left.
{}-
  \left(
    \frac{ g(x, u  (x)) }{ u  (x)^{r-1} }
  - \frac{ g(x, u_0(x)) }{ u_0(x)^{r-1} }
  \right)
  \right]
  (u(x)^r - u_0(x)^r) \,\mathrm{d}x\leq 0 \,,
\end{aligned}
\end{align*}
as the function
\begin{math}
  s\mapsto [ f(x,s) - g(x,s) ] / s^{r-1} \colon (0,\infty)\to \RR_+
\end{math}
is strictly monotone decreasing for every $x\in \Omega$,
by Hypotheses {\bf (f2)} and~{\bf (g2)}.
Since the opposite inequality ``$\geq$'' must be valid,
by ineq.~\eqref{in_A:Diaz-Saa},
we conclude that the equality above must hold.
This forces
$u(x) = u_0(x)$ at almost every point $x\in \Omega$,
by Hypotheses {\bf (f2)} and~{\bf (g2)}, i.e.,
$u\equiv u_0$ in $\Omega$.

The proof of Theorem~\ref{thm-problem_2} is finished.
\qed
\par\vskip 10pt

Theorem~\ref{thm-problem_2} has the following interesting special case.

\begin{cor}\label{cor-problem_2}
Assume that\/
$p\in C^{0,\alpha_1}(\overline{\Omega})$ for some $\alpha_1\in (0,1)$
and the constant $r\in [1,\infty)$ satisfy
{\rm Hypothesis ($\mathbf{p}$)} together with $r < p_{-}$, i.e.,
\begin{equation*}
  1\leq r < p_{-}\eqdef \inf_{\Omega} p(x)\leq
            p_{+}\eqdef \sup_{\Omega} p(x) < \infty \,.
\end{equation*}
Let\/
$h, \ell\in C(\overline{\Omega})$ and\/
$q, Q\in C(\overline{\Omega})$
be two pairs of strictly positive functions such that\/
\begin{align}
\label{e:q_+<r<p-}
&
\left\{
\qquad
\begin{aligned}
  1\leq {}
&      q_{-}\eqdef \inf_{\Omega} q(x)\leq
       q_{+}\eqdef \sup_{\Omega} q(x)
\\
  < r < {}
&      p_{-}\eqdef \inf_{\Omega} p(x)\leq
       p_{+}\eqdef \sup_{\Omega} p(x) < \infty \,,
\end{aligned}
\right.
\\
\label{e:r<Q_}
&      r\leq Q_{-}\eqdef \inf_{\Omega} Q(x)\leq
             Q_{+}\eqdef \sup_{\Omega} Q(x) < \infty \,.
\end{align}
Let\/
$f,g\colon \overline{\Omega}\times \RR_+\to \RR_+$
be defined by\/
$f(x,s) = h(x)\, s^{q(x) - 1}$ and\/
$g(x,s) = \ell(x)\, s^{Q(x) - 1}$ for\/
$(x,s)\in \overline{\Omega}\times \RR_+$.
Then the conclusion of {\rm Theorem~\ref{thm-problem_2}}
for problem~\eqref{e:problem_2}
taking the following special form,
\begin{equation}
\label{eq:problem_2}
\left\{
\begin{alignedat}{2}
  {}- \Delta_{p(x)} u + \ell(x)\, u^{Q(x) - 1}
& {}= h(x)\, u^{q(x) - 1}
&&  \quad\mbox{ in }\, \Omega \,;
\\
  u & {}= 0
&&  \quad\mbox{ on }\, \partial\Omega \,,
    \quad u > 0 \mbox{ in }\, \Omega \,,
\end{alignedat}
\right.
\end{equation}
\end{cor}

\par\vskip 10pt
\proof
It is a matter of easy, direct calculations that functions
$f$ and $g$ satisfy all Hypotheses {\bf (f1)} -- {\bf (f3)}
and all Hypotheses {\bf (g1)} -- {\bf (g3)}, respectively.
Notice that
\begin{math}
  1\leq q(x) < r\leq \min\{ p(x) ,\, Q(x)\}
\end{math}
holds for all $x\in \overline{\Omega}$.
\qed
\par\vskip 10pt

Our last application concerns a nonlocal boundary value problem
of {\it Kirchhoff's type\/} involving local and nonlocal nonlinearities
treated e.g.\ in
{\sc Ch.-Y.\ Chen}, {\sc Y.-Ch.\ Kuo}, and {\sc Ts.-F.\ Wu}
\cite{Chen-Kuo-Wu}.
This problem is motivated by the stationary (elliptic) case of
an evolutionary hyperbolic equation that arises in the study of
string or membrane vibrations, where $u = u(x,t)$ stands for
the displacement at $x\in \Omega$ and time $t\in \RR_+$, cf.\
\cite[Eq.\ (1.2), p.~1877]{Chen-Kuo-Wu}.
The mathematical model for the stationary displacement
$u = u(x)$ at $x\in \Omega$ takes the following form,
\begin{equation}
\label{e:problem_3}
\left\{
\begin{alignedat}{2}
  {}- M
    \left( \int_{\Omega}
  \frac{ \left\vert \nabla u(x) \right\vert^{p(x)} }{p(x)}\,
  \,\mathrm{d}x
    \right)\, \Delta_{p(x)} u
& {}= f(x,u)
&&  \quad\mbox{ in }\, \Omega \,;
\\
  u & {}= 0
&&  \quad\mbox{ on }\, \partial\Omega \,,
    \quad u > 0 \mbox{ in }\, \Omega \,.
\end{alignedat}
\right.
\end{equation}
In the original physics problem, $p(x)\equiv 2$ is constant.

\par\vskip 10pt
%
In addition to our Hypotheses {\bf (f1)} -- {\bf (f3)} imposed on
the function $f$ at the beginning of this section,
we impose the following hypotheses on the function $M$:
\begin{itemize}
\item[{\bf (M1)}]
$\;$
$M\colon \RR_+\to \RR_+$
is a nonnegative continuous function with $M(0) > 0$.
\item[{\bf (M2)}]
$\;$
$M\colon \RR_+\to \RR_+$ is monotone increasing,
but {\em not\/} necessarily {\em strictly\/} monotone increasing.
\item[{\bf (M3)}]
$\;$
$M\colon \RR_+\to \RR_+$ is bounded, that is, the monotone limit
$M(s)\nearrow M(+\infty) < \infty$ as $s\nearrow +\infty$
is finite.
\end{itemize}
%
\par\vskip 10pt

As a consequence of Hypothesis~{\bf (M2)} we obtain also
\begin{equation}
\label{e:M_hat}
  M(0)\, t\leq
  \hat{M}(t)\eqdef \int_0^t M(s) \mathrm{d}s
  \leq M(+\infty)\, t
  \quad\mbox{ for every }\, t\in \RR_+ \,.
\end{equation}
Clearly,
$\hat{M}\colon \RR_+\to \RR_+$
is strictly monotone increasing and convex
(possibly {\em not\/} strictly convex).
Recalling the potential $F$ introduced in eq.~\eqref{e:F'=f-probl_1},
we observe that problem \eqref{e:problem_3}
corresponds to the Euler equation for a critical point
$u\in W_0^{1,p(x)}(\Omega)$ of the energy functional
$\mathcal{J}\colon W_0^{1,p(x)}(\Omega)\to \RR_+$ defined by
\begin{equation}
\label{def:J(u)}
  \mathcal{J}(u)\equiv \mathcal{J}_{p(x),f}(u)\eqdef
  \hat{M}
    \left( \int_{\Omega}
  \frac{ \left\vert \nabla u(x) \right\vert^{p(x)} }{p(x)}\,
  \,\mathrm{d}x
    \right)
  - \int_{\Omega} F(x,u(x)) \,\mathrm{d}x
\end{equation}
for every function $u\in W_0^{1,p(x)}(\Omega)$.
This functional is well\--defined, by the Sobolev embedding
$W_0^{1,p(x)}(\Omega) \hookrightarrow L^r(\Omega)$,
which is even compact,
and by the estimate in \eqref{lim:s=infty}.
The reader is referred to the monograph by
{\sc L.\ Diening}, {\sc P.\ Harjulehto}, {\sc P.\ H\"ast\"o}, and
{\sc M.\ R\r{u}\v{z}i\v{c}ka}
\cite[{\S}8.3 and {\S}8.4]{DienHHR}
for Sobolev embeddings and their compactness.
It is coercive thanks to inequalities
\eqref{lim:s=infty}, \eqref{e:M_hat}, and $r\leq p_{-}$, i.e.,
it satisfies an analogue of \eqref{e:E-coerce},
\begin{equation*}
  \| u\|_{ W^{1,p(x)}_0(\Omega) }\eqdef
  \|\nabla u\|_{ L^{p(x)}(\Omega) } \,\longrightarrow\, +\infty
  \quad\Longrightarrow\quad \mathcal{J}(u)\to +\infty \,.
\end{equation*}
Furthermore, it is easy to see that
$\mathcal{J}$ is G\^ateaux\--differentiable on $W_0^{1,p(x)}(\Omega)$
with the G\^ateaux derivative
\begin{align}
\label{def:J'(u)}
&
\begin{aligned}
  \left[ \mathcal{J}'(u)\right] (x) = {}
& \hat{M}'
    \left( \int_{\Omega}
  \frac{ \left\vert \nabla u(x) \right\vert^{p(x)} }{p(x)}\,
  \,\mathrm{d}x
    \right)\cdot
  \left[ {}-
  \mathrm{div} \left( |\nabla u|^{p(x)-2} \nabla u\right)
  \right]
\\
& {}
  - f(x,u)
\end{aligned}
\\
\nonumber
& {}
  = M
    \left( \int_{\Omega}
  \frac{ \left\vert \nabla u(x) \right\vert^{p(x)} }{p(x)}\,
  \,\mathrm{d}x
    \right)\cdot (- \Delta_{p(x)} u) - f(x,u) \,,\quad x\in \Omega \,,
\end{align}
which may be interpreted as a distribution over $\Omega$,
that is, it belongs to the locally convex space
$\mathcal{D}'(\Omega)$ of all distributions over $\Omega$
which is the dual space of
\begin{math}
  \mathcal{D}(\Omega) = C^{\infty}_{\mathrm{c}}(\Omega) .
\end{math}

We have the following analogue of
Theorems \ref{thm-problem_1} and~\ref{thm-problem_2}
for positive weak solutions
to the boundary value problem~\eqref{e:problem_3}:

\begin{theorem}\label{thm-problem_3}
Under the {\rm Hypotheses\/}
{\rm ($\boldsymbol{\Omega}$)}, {\rm ($\mathbf{p}$)},
{\bf (f1)} -- {\bf (f3)}, and\/ {\bf (M1)} -- {\bf (M3)},
the nonlocal Kirchhoff problem~\eqref{e:problem_3} possesses
a unique nonnegative and nontrivial weak solution
$u\in W^{1,p(x)}_0(\Omega)\cap L^\infty(\Omega)$.
This solution belongs to the class
$C^{1,\beta}(\overline{\Omega})$, for some $\beta\in (0, \alpha)$,
and satisfies also the Hopf maximum principle
\eqref{Hopf:problem_1},
\begin{equation*}
  u(x) > 0 \;\mbox{ for all }\, x\in \Omega
  \quad\mbox{ and }\quad
  \frac{\partial u}{\partial\boldsymbol{\nu}} (x) < 0
    \;\mbox{ for all }\, x\in \partial\Omega \,.
\end{equation*}
Of course, $u = 0$ on the boundary $\partial\Omega$.
Hence, $u$ is also a positive weak solution.
\end{theorem}

\par\vskip 10pt
\proof
Although we could generalize
the {\sc D{\'\i}az} and {\sc Saa} inequality \eqref{ineq:Diaz-Saa}
(proved in Theorem~\ref{thm-Diaz-Saa})
to the class of nonlocal quasilinear elliptic operators
as suggested in the Kirchhoff problem~\eqref{e:problem_3},
we prefer to give a direct proof of our theorem
which, however, follows very closely the same ideas as does
our proof of Theorem~\ref{thm-Diaz-Saa}.

We begin with the following trivial observation;
we use the same notation as does our convexity result in
Theorem~\ref{thm-Diaz-Saa}:

The composition functional
$\hat{M}\circ \mathcal{W}_{A}\colon W\to \RR_+$ is given by
\begin{align}
\label{def:M.W_A(v)}
& \left[ \hat{M}\circ \mathcal{W}_A\right] (v)
  \equiv \hat{M}\left( \mathcal{W}_{A,p(x),r}(v) \right)
\\
\nonumber
& {}
  = \hat{M}\left(
    \int_{\Omega} \frac{r}{p(x)}\,
    A \left( x,\, \frac{ \nabla (|v|^{1/r}) }%
                  { \left\vert \nabla (|v|^{1/r}) \right\vert }
      \right)\cdot
    \left\vert \nabla (|v|^{1/r}) \right\vert^{p(x)} \,\mathrm{d}x
  \right)
\end{align}
for every function $v\in W$; see eqs.\
\eqref{def:W(v)} (in the Introduction, Section~\ref {s:Intro}),
\eqref{def:W_A(v)}, and~\eqref{e:M_hat}.
In particular, concerning the Kirchhoff problem~\eqref{e:problem_3},
we take
$A(x,\xi) = |\xi|^{p(x)}$ for $(x,\xi)\in \Omega\times \RR^N$,
in which case $A(x,\xi) = 1$ for all
$(x,\xi)\in \Omega\times \SS^{N-1}$.

By our convexity result in Theorem~\ref{thm-ray-convex},
the restriction of the functional\/
$\mathcal{W}_A\colon W\to \RR_+$ to the convex cone $\dotV$ is
ray\--strictly convex on~$\dotV$.
Recall from above that
$\hat{M}\colon \RR_+\to \RR_+$
is strictly monotone increasing and convex.
Consequently, an easy exercise in elementary analysis reveals
that also the composition functional
$\hat{M}\circ \mathcal{W}_{A}\colon W\to \RR_+$
must be ray\--strictly convex on~$\dotV$.
By our Hypotheses {\bf (f1)} and {\bf (f2)} on~$f$,
for every fixed $x\in \Omega$, also the function
$t\mapsto {}- F(x, t^{1/r})\colon \RR_+\to \RR_+$
is strictly convex, owing to the partial derivative
\begin{equation*}
    t\mapsto
    \frac{\partial}{\partial t}\, F(x, t^{1/r})
  = \frac{1}{r}\, t^{-1 + (1/r)}\cdot
    \frac{\partial F}{\partial s} (x, t^{1/r})
  = \frac{1}{r}\cdot
    \frac{ f(x, t^{1/r}) }{ (t^{1/r})^{r-1} }
  \colon (0,\infty)\to \RR_+
\end{equation*}
being strictly monotone decreasing on $(0,\infty)$.
From these two convexity results we deduce that also the functional
\begin{align}
\label{def:M.W_A(v)+F}
  v \;\longmapsto\;
& \hat{\mathcal{J}}(v)\eqdef
    \mathcal{J}\left( |v|^{1/r} \right)
  = \left[ \hat{M}\circ \mathcal{W}_A\right] (|v|^{1/r})
  - \int_{\Omega} F\left( x, |v(x)|^{1/r} \right) \,\mathrm{d}x
\\
\nonumber
& {}
  = \hat{M}\left(
    \int_{\Omega} \frac{r}{p(x)}\cdot
    \left\vert \nabla (|v|^{1/r}) \right\vert^{p(x)} \,\mathrm{d}x
  \right)
  - \int_{\Omega} F\left( x, |v(x)|^{1/r} \right) \,\mathrm{d}x
\end{align}
must be strictly convex on~$\dotV$.

Now we are ready to prove the uniqueness claim in our theorem:
On the contrary, let us assume that
$u_1, u_2\in W^{1,p(x)}_0(\Omega)\cap L^\infty(\Omega)$
are two different nonnegative and nontrivial weak solutions
to problem~\eqref{e:problem_3} that satisfy also
the Hopf maximum principle \eqref{Hopf:problem_1}.
In particular, we have
$\mathcal{J}'(u_1) = \mathcal{J}'(u_2) = 0$
in $\mathcal{D}'(\Omega)$.
Setting $v_1 = u_1^r$ and $v_2 = u_2^r$ we get also
the G\^ateaux derivatives
$\hat{\mathcal{J}}'(v_1) = \hat{\mathcal{J}}'(v_2) = 0$
as distributions in $\mathcal{D}'(\Omega)$.
Moreover, we have $v_1, v_2\in \dotV$ and
$v_1/v_2$, $v_2/v_1\in L^{\infty}(\Omega)$.
Consequently, also
$v\eqdef (1-\theta) v_1 + \theta v_2\in \dotV$
is valid for all $\theta\in (-\delta, 1+\delta)$,
where $\delta\in (0,1)$ is small enough.
The function
\begin{equation*}
  \theta\mapsto \Phi(\theta)\eqdef \hat{\mathcal{J}}(v)
  = \hat{\mathcal{J}}\left( (1-\theta) v_1 + \theta v_2\right)
    \colon (-\delta, 1+\delta)\to \RR_+
\end{equation*}
is strictly convex and differentiable with the derivative
\begin{align}
\label{e:dJ/d_theta}
\begin{aligned}
  \Phi'(\theta)
& {}
  = M\left(
    \int_{\Omega} \frac{r}{p(x)}\cdot
    \left\vert \nabla (|v|^{1/r}) \right\vert^{p(x)} \,\mathrm{d}x
  \right)
\\
& {}
  \times \int_{\Omega}
  \left\vert \nabla (|v(x)|^{1/r}) \right\vert^{p(x) - 2}
  \nabla (|v(x)|^{1/r})\cdot
  \nabla \genfrac{(}{)}{}0{v_2-v_1}{ v^{ 1 - \frac{1}{r} } }
  \,\mathrm{d}x
\end{aligned}
\\
\nonumber
  {}
  - \frac{1}{r}
    \int_{\Omega} f(x, |v(x)|^{1/r})\cdot
  \genfrac{}{}{}0{v_2-v_1}{ v^{ 1 - \frac{1}{r} } }
  \,\mathrm{d}x \,.
\end{align}
The monotonicity of the derivative
\begin{math}
  \theta\mapsto \Phi'(\theta)\colon (-\delta, 1+\delta)\to \RR
\end{math}
yields
\begin{equation*}
  0\leq \Phi'(t) - \Phi'(0)\leq \Phi'(1) - \Phi'(0)
  \quad\mbox{ for every }\, t\in [0,1] \,.
\end{equation*}
But
$\hat{\mathcal{J}}'(v_1) = \hat{\mathcal{J}}'(v_2) = 0$
in $\mathcal{D}'(\Omega)$ forces $\Phi'(0) = \Phi'(1) = 0$
whence $\Phi'(t) = 0$ for every $t\in [0,1]$.
We conclude that
$\Phi(t) = \Phi(0)$ for every $t\in [0,1]$
which contradicts the strict convexity of $\Phi$ on $[0,1]$.

The uniqueness part of our theorem is proved.

To verify the existence part,
we apply analogous arguments as in our proofs of Theorems
\ref{thm-problem_1} and \ref{thm-problem_2}.
Recalling that the energy functional $\mathcal{J}$ defined in
eq.~\eqref{def:J(u)} is coercive and weakly lower semi\-continuous on
$W^{1,p(x)}_0(\Omega)$, by
\cite[{\S}13.2, pp.\ 412--417]{DienHHR},
we conclude that it possesses
a global minimizer $u_0\in W^{1,p(x)}_0(\Omega)$, by
\cite[Theorem 1.2, p.~4]{Struwe}.
To verify $u_0\geq 0$ a.e.\ in $\Omega$,
we first observe that also its positive part,
$u_0^{+}\eqdef \max\{ u_0, 0\} \geq 0$,
is in $W_0^{1,p(x)}(\Omega)$ and thus satisfies
$\mathcal{J}(u_0^{+}) \geq \mathcal{J}(u_0)$.
Recalling
$\Omega^{+} = \{ x\in \Omega\colon u_0(x)\geq 0\}$ and
$\Omega^{-} = \{ x\in \Omega\colon u_0(x)\leq 0\}$, we calculate
\begin{align*}
&
\begin{aligned}
    \mathcal{J}(u_0)
& {}
  = \hat{M} \left( \int_{\Omega}
  \frac{ \left\vert \nabla u_0(x) \right\vert^{p(x)} }{p(x)}\,
  \,\mathrm{d}x
    \right)
  - \int_{\Omega} F(x,u_0(x)) \,\mathrm{d}x
\end{aligned}
\\
&
\begin{aligned}
{}= \hat{M}
    \left( \int_{ \Omega^{+} }
  \frac{ \left\vert \nabla u_0(x) \right\vert^{p(x)} }{p(x)}\,
  \,\mathrm{d}x
         + \int_{ \Omega^{-} }
  \frac{ \left\vert \nabla u_0(x) \right\vert^{p(x)} }{p(x)}\,
  \,\mathrm{d}x
    \right)
\\
{}- \int_{ \Omega^{+} } F(x,u_0(x)) \,\mathrm{d}x
\end{aligned}
\\
&
\begin{aligned}
& {}
  \geq \hat{M}
    \left( \int_{ \Omega^{+} }
  \frac{ \left\vert \nabla u_0(x) \right\vert^{p(x)} }{p(x)}\,
  \,\mathrm{d}x
    \right)
  - \int_{ \Omega^{+} } F(x,u_0(x)) \,\mathrm{d}x
\end{aligned}
\\
&
\begin{aligned}
& {}
  = \mathcal{J}(u_0^{+}) \geq \mathcal{J}(u_0) \,.
\end{aligned}
\end{align*}
In fact, these inequalities must be equalities.
Since $\hat{M}\colon \RR_+\to \RR_+$
is strictly monotone increasing with $\hat{M}'= M >0$ in $\RR_+$,
by {\rm Hypothesis\/} {\bf (M1)},
the equalities above force
$\nabla u_0^{-}(x) = {}- \nabla u_0(x)\equiv 0$
a.e.\ in $\Omega^{-}$, whence
$u_0^{-}(x)\equiv 0$ a.e.\ in~$\Omega^{-}$.
We have proved $u_0\geq 0$ a.e.\ in $\Omega$ as claimed.

In order to exclude the possibility that $u_0\equiv 0$ in $\Omega$,
we construct a function
$u_1\in W_0^{1,p(x)}(\Omega)$ such that
$\mathcal{J}(u_1) < 0 = \mathcal{J}(0)$.
To this end, we take an arbitrary nonnegative $C^1$-function
$\phi\in C^1_{\mathrm{c}}(\Omega)$
with compact support in $\Omega$, $\phi\not\equiv 0$ in $\Omega$.
In analogy with ineq.~\eqref{e_hat:E-zero},
for $0 < t\leq 1$ we invoke
{\rm Hypothesis\/} {\bf (M1)} to estimate
\begin{align}
\label{e:J-zero}
    \mathcal{J}(t\phi)
& {}
  = \hat{M}
    \left( \int_{\Omega} \frac{ t^{p(x)} }{p(x)}\,
    |\nabla\phi(x)|^{p(x)} \,\mathrm{d}x
    \right)
  - \int_{\Omega} F(x, t\phi(x)) \,\mathrm{d}x
\\
\nonumber
& \leq \hat{M}
    \left( \frac{ t^{p_{-}} }{p_{-}}
  \int_{\Omega} |\nabla\phi(x)|^{p(x)} \,\mathrm{d}x
  \right)
  - \int_{\Omega} F(x, t\phi(x)) \,\mathrm{d}x \,.
\end{align}
Recall that
$1\leq r\leq p_{-} = \inf_{\Omega} p(x)$,
by our {\rm Hypothesis ($\mathbf{p}$)}, ineq.~\eqref{e:p-<p<p+}.
We take advantage of inequalities
\eqref{e:M_hat} (for $\hat{M}$) and \eqref{lim_F:s=0} (for $F$)
to estimate the last two terms in ineq.~\eqref{e:J-zero} above,
\begin{align*}
    \mathcal{J}(t\phi)
& {}
  \leq M(+\infty)\, \frac{t^r}{r}
  \int_{\Omega} |\nabla\phi(x)|^{p(x)} \,\mathrm{d}x
  - \frac{t^r}{r\eps}
  \int_{\Omega} \phi(x)^r \,\mathrm{d}x
\\
& = {}- \frac{t^r}{r}
    \left(
    \frac{1}{\eps} \int_{\Omega} \phi(x)^r \,\mathrm{d}x
  - M(+\infty) \int_{\Omega} |\nabla\phi(x)|^{p(x)} \,\mathrm{d}x
    \right)
\end{align*}
for all $t\in [0,t_{\eps}]$.
Choosing $\eps > 0$ small enough, we conclude that
$\mathcal{J}(t\phi) < 0$ whenever $0 < t\leq t_{\eps}$.
In addition to $u_0\geq 0$ a.e.\ in $\Omega$,
we have proved also $u_0\not\equiv 0$ in~$\Omega$.

Since $u_0\in W_0^{1,p(x)}(\Omega)$
is a global minimizer for the functional
$\mathcal{J}\colon W_0^{1,p(x)}(\Omega)\to \RR$,
it is also a critical point for $\mathcal{J}$ and, hence,
a nonnegative weak solution to problem~\eqref{e:problem_3}
provided $u_0\in L^{\infty}(\Omega)$.
As
\begin{equation*}
  M_0 = M\left( \int_{\Omega}
  \frac{ \left\vert \nabla u_0(x) \right\vert^{p(x)} }{p(x)}\,
  \,\mathrm{d}x
    \right)
\end{equation*}
is a positive constant,
$0 < M(0)\leq M_0\leq M(+\infty) < \infty$,
the Dirichlet problem~\eqref{e:problem_3} for $u = u_0$
is identical with that in~\eqref{e:problem_1}
with $f(x,u_0)$ replaced by $M_0^{-1}\, f(x,u_0)$.
The rest of the proof now follows from the conclusion of
Theorem~\ref{thm-problem_1}.
\qed
\par\vskip 10pt

\begin{remark}\label{rem-Motreanu}\nopagebreak
\begingroup\rm
The question of uniqueness of general,
possibly sign\--changing weak solutions to problems of type
\eqref{e:problem_1} was studied in a number of articles; see, e.g.,
{\sc S.~N.\ Antontsev}, {\sc M.\ Chipot}, and {\sc Y.\ Xie}
\cite{Ant-Chipot},
{\sc S.~N.\ Antontsev} and {\sc S.~I.\ Shmarev} \cite{Ant-Shm-1},
{\sc V.~V.\ Motreanu} \cite{Motreanu}, and references therein.
There, the function
$f(x,u)$ on the right\--hand side of our problem~\eqref{e:problem_1}
is replaced by a somewhat more general reaction function
$f(x) - b(x,u)$ for $x\in \Omega$ and $u\in \RR$, where
$s\mapsto b(x,s)\colon \RR\to \RR$
is assumed to be continuous and monotone increasing
(i.e., nondecreasing) for almost every $x\in \Omega$.
Some of the most recent results in \cite{Motreanu} require that
$s\mapsto b(x,s)$ be even strictly monotone increasing, although
two weak comparison results in \cite[Theorems 4.1 and 4.2]{Motreanu}
are proved for $s\mapsto b(x,s)$ being nondecreasing only.

All these results are based on the fact that the quasilinear operator
$u\mapsto {}- \Delta_{p(x)} u + b(x,u(x))$
is {\em monotone\/} with respect to the $L^2(\Omega)$\--induced duality
between $W_0^{1,p(x)}(\Omega)$ and its dual space.
Of course, a mild additional condition guaranteeing some kind of
{\em strict monotonicity\/} must be imposed in order to obtain
the desired uniqueness result.
These hypotheses on the quasilinear operator may be reformulated
in terms of convexity properties of the corresponding energy functional;
cf.\ eq.~\eqref{def:E(v)}.

According to the original observation in
{\sc H.\ Br\'ezis} and {\sc L.\ Oswald} \cite{Brezis-Osw}
for $p(x)\equiv 2$ (a constant),
much stronger convexity properties of this energy functional
can be proved for {\em positive solutions\/}
to problem~\eqref{e:problem_1}:
Namely, those first observed in
{\sc J. Fleckinger} et al.\ \cite{FleckHTT}
for $p(x)\equiv p\in (1,\infty)$ (a constant) and generalized
in the present article for a variable exponent $p(x)$.
The mechanism of this approach, based on \cite{Brezis-Osw},
composes the standard convex energy functional
\begin{math}
  \mathcal{W}_{p(x),1}\colon W_0^{1,p(x)}(\Omega)\to \RR_+
\end{math}
defined in eq.~\eqref{def:W(v)} for $r=1$
with the Nemytskii operator given by the concave function
$s\mapsto s^{1/r}\colon \RR_+\to \RR_+$
($1 < r\leq p_{-} = \inf_{\Omega} p$).
Somewhat surprisingly, when restricted to
the (convex) cone of positive functions in $W_0^{1,p(x)}(\Omega)$,
this composition is still convex, even {\em ray\--strictly convex},
as proved in our Theorem~\ref{thm-ray-convex}.
Unlike in {\sc V.~V.\ Motreanu} \cite{Motreanu},
we are then able to treat problem~\eqref{e:problem_1} with 
$s\mapsto b(x,s)\colon \RR\to \RR$
strictly monotone decreasing (in the notation of \cite{Motreanu});
see Hypothesis~{\bf (f2)} and Example~\ref{exam-appl_1}.
\endgroup
\end{remark}
\par\vskip 10pt


\section*{Acknowledgement.}
The work of Peter Tak\'a\v{c}
was partially supported by a grant from
Deutsche Forschungs\-gemeinschaft (DFG, Germany),
Grant no.\ TA~213/16-1.
The work of Jacques Giacomoni was supported in part by
C.N.R.S. (France) through
``Laboratoire de Math\'ematiques et de leurs Applications -- Pau'',
no.\ $5142$.
Both authors would like to express their thanks to
an anonymous referee for pointing out several unclear formulations
and giving us hints to new pertinent references.
His/her careful and detailed comments have helped to improve
the manuscript in a significant way.



%
%
\makeatletter \renewcommand{\@biblabel}[1]{\hfill#1.} \makeatother
%
%

\end{document}